\documentstyle{article}

\newcommand{\ms}{\medskip}

\newcommand{\noi}{\noindent}
\newcommand{\ra}{\rightarrow}
\newcommand{\bea}{\begin{eqnarray}}
\newcommand{\eea}{\end{eqnarray}}
\newcommand{\ol}{\overline}
\newcommand{\gr}{Groenewold}
\newcommand{\vh}{Van~Hove}
\newcommand{\vn}{Von~Neumann}
\newcommand{\q}{{\cal Q}}
\newcommand{\p}{{\cal P}}
\newcommand{\h}{{\cal H}}
\newcommand{\f}{{\cal F}}

\newcommand{\A}{{\cal A}}
\newcommand{\bb}{{\cal B}}
\newcommand{\oo}{{\cal O}}
\newcommand{\jj}{{\cal J}}
\newcommand{\s}{{\cal S}}
\newcommand{\uu}{{\cal U}}
\newcommand{\nn}{{\cal N}}

\def\r{{\bf R}}
\def\z{{\bf Z}}

\newtheorem{thm}{Theorem}
\newtheorem{lem}{Lemma}
\newtheorem{cor}{Corollary}
\newtheorem{prop}[thm]{Proposition}
\newtheorem{defn}{Definition}
\newtheorem{conj}{Conjecture}

\def\f #1,#2.{\textstyle{#1\over #2}}
\def\hlf{{\f 1,2.}}
\def\C{{\bf C}}
\def\sp{{\rm span}}

\title{{Obstruction Results in Quantization Theory}\thanks{To appear in J.
Nonlinear Sci.}}

\author{{\bf Mark J. Gotay}\thanks{Supported in part by NSF grants DMS
92-22241 and 96-23083.}
\\  Department of Mathematics \\ University of Hawai`i \\ 2565 The
Mall \\ Honolulu, HI 96822  USA \\
\and {\bf Hendrik B. Grundling} \\ Department of Pure Mathematics \\
University of New South Wales \\ P.O. Box 1 \\ Kensington, NSW 2033
Australia \\
\and  {\bf Gijs M. Tuynman} \\ URA 751 au CNRS \& UFR de
Math\'ematiques \\ Universit\'{e} de Lille I \\ F-59655 Villeneuve
d'Ascq cedex France}

\date{\today}

\begin{document}


\maketitle

\begin{abstract} Quantization is not a straightforward proposition,
as demonstrated by Groenewold's and Van~Hove's discovery, exactly
fifty years ago, of an ``obstruction'' to quantization. Their ``no-go
theorems'' assert that it is in principle  impossible to consistently
quantize every classical observable on the phase space ${\bf R}^{2n}$
in a physically meaningful way. A similar obstruction was recently
found for
$S^2$, buttressing the common belief that no-go theorems should hold
in some generality. Surprisingly, this is not so -- it has also just
been proven that there is no obstruction to quantizing a torus.

In this paper we take first steps towards delineating the
circumstances under which such obstructions will appear, and
understanding the mechanisms which produce them. Our objectives are
to conjecture a generalized \gr-\vh\ theorem, and to determine the
maximal subalgebras of observables which can be consistently
quantized. This requires a study of the structure of Poisson algebras
of classical systems and their representations. To these ends we
include an exposition of both prequantization (in an extended sense)
and quantization theory -- formulated in terms of ``basic sets of
observables,'' and review in detail the known results for
$\r^{2n}$, $S^2$, and $T^2$. Our discussion is independent of any
particular method of quantization; we concentrate on the structural
aspects of quantization theory which are common to all Hilbert
space-based quantization techniques.

\end{abstract}


\begin{section}{Introduction}

Quantization -- the problem of  constructing the quantum formulation
of a system from its classical description -- has always been one of
the great mysteries of mathematical physics. It is generally
acknowledged that quantization is an ill-defined procedure, which
cannot be consistently applied to all classical systems. While there
is certainly no extant quantization procedure which works well in all
circumstances, this assertion nonetheless bears closer scrutiny.

Already from first principles one encounters difficulties. Given that
the  classical description of a system is an approximation to its
quantum description, obtained in a macroscopic limit (when
$\hbar\to 0$), one expects that some information is lost in the
limit. So quantization should somehow have to compensate for this.
But how can a given quantization procedure select, from amongst the
myriad of quantum theories all of which have the same classical
limit, the physically correct one?

In view of this ambiguity it is not surprising that the many
quantization schemes which have been developed over the years -- such
as the physicists' original ``canonical quantization''
\cite{di} (and its modern formulations, such as geometric quantization
\cite{ki,So70,Wo}), Weyl quantization \cite{Fo} (and its successor
deformation quantization
\cite{Bayen e.a.,ri2,ri3}), path integral quantization \cite{GJ}, and
the group theoretic approach to quantization \cite{Is}, to cite just
some -- have shortcomings. Rather, is it amazing that they work as
well as they do!

But there are deeper, subtler problems, involving the Poisson
algebras of classical systems and their representations. In this
context the conventional wisdom is that it is impossible to ``fully''
quantize any given classical system -- regardless of the particular
method employed -- in a way which is consistent with the physicists'
Schr\"odinger quantization of
${\bf R}^{2n}$. (We will make this somewhat nebulous statement precise
later.) In other words, the assertion is that there exists a universal
``obstruction'' which forces one to settle for something less than a
complete and consistent quantization of {\it any\/} system. Each
quantization procedure listed above evinces this defect in various
examples.

That there are problems in quantizing even simple systems was
observed very early on. One difficulty was to identify the analogue
of the multiplicative structure of the classical observables in the
quantum formalism. For instance, consider the quantization of ${\bf
R}^{2n}$ with canonical coordinates
$\{q^i,\,p_i\,|\,i = 1,\dots,n\}$, representing the phase space of a
particle moving in $\r^n$. For simple observables the ``product
$\ra$ anti-commutator'' rule worked well. But for more complicated
observables (say, ones which are quartic polynomials in the positions
and momenta), this rule leads to inconsistencies. (See
\cite[\S 4]{a-b},
\cite[\S1.1]{Fo} and \S 4 for a discussion of these factor-ordering
ambiguities.) Of course this, in and by itself, might only indicate
the necessity of coming up with some subtler symmetrization rule. But
attempts to construct a quantization map also conflicted with Dirac's
``Poisson bracket $\ra$ commutator'' rule. This was implicitly
acknowledged by Dirac
\cite[p. 87]{di}, where he made the now famous hedge:

\begin{quote}
\em ``The strong analogy between the quantum P.B.
$\ldots$ and the classical P.B.
$\ldots$ leads us to make the assumption that the quantum P.B.s, {\rm
or at any rate the simpler ones of them,} have the same values as the
corresponding classical P.B.s.''
\end{quote}

\noi In any case, as a practical matter, one was forced to limit the
quantization to relatively ``small'' subalgebras of classical
observables which could be handled without ambiguity (e.g.,
polynomials which are at most quadratic in the $p$'s and the
$q$'s, or observables which are at most affine functions of the
coordinates or of the momenta).

Then, in 1946, Groenewold \cite{Gr} showed that the search for an
``acceptable'' quantization map was futile. The strong version of his
``no-go'' theorem states that one cannot consistently quantize the
Poisson algebra of all polynomials in the $q^i$ and $p_i$ on ${\bf
R}^{2n}$ as symmetric operators on some Hilbert space $\h,$ subject to
the requirement that the
$q^i$ and  $p_i$ be irreducibly represented.\footnote{\,There are
actually two variants of \gr's theorem (``strong'' and ``weak''); both
will be discussed in \S 4.1.} Van~Hove subsequently refined \gr's
result \cite{vH1}. Thus it is {\it in principle} impossible to
quantize -- by {\it any\/} means -- every classical observable on
${\bf R}^{2n}$, or even every polynomial observable, in a way
consistent with Schr\"odinger quantization (which, according to the
Stone-\vn\ theorem, is the import of the irreducibility requirement
on the $p$'s and $q$'s). At most, one can consistently quantize
certain subalgebras of observables, for instance the ones mentioned
in the preceding paragraph.


Of course, \gr's remarkable result is valid only for the classical
phase space
$\r^{2n}.$ The immediate problem is to determine whether similar
obstructions appear when trying to quantize other symplectic
manifolds. Little is known in this regard, and only in the past few
years have other examples come to light. Just recently an obstruction
was found in the case of the symplectic manifold $S^2$, representing
the (internal) phase space of a massive spinning particle \cite{GGH}.
It was shown that one cannot consistently quantize the Poisson algebra
of spherical harmonics (thought of as ``polynomials'' in the components
$S_i$ of the spin angular momentum vector $\bf S$), subject to the
requirement that the
$S_i$ be irreducibly represented on a Hilbert space of dimension
greater than one. This is a direct analogue for $S^2$ of
\gr's theorem. Combined with the observation that
$S^2$ is in a sense at the opposite extreme from
$\r^{2n}$ insofar as symplectic manifolds go, it indicates that no-go
theorems can be expected to hold in some generality. But,
interestingly enough, they are {\em not\/} universal: it is possible
to explicitly exhibit a quantization of the full Poisson algebra of
the torus $T^2$ in which a suitable irreducibility requirement is
imposed \cite{Go}. An important point, therefore, is to understand
the mechanisms which are responsible for these divergent outcomes.

Our goal here is to study such obstruction results for the
quantization of a Poisson algebra of a symplectic manifold. We will
review the known results in some detail, and give a careful
presentation of prequantization (in an extended sense) and
quantization, with a view to conjecturing a generalized \gr-\vh\
theorem and in particular delineating the circumstances under which
it can be expected to hold. Our discussion will be independent of any
particular method of quantization; we concentrate on the structural
aspects of quantization theory which are common to all Hilbert
space-based quantization techniques.

\end{section}


\begin{section}{Prequantization}

Let
$(M,\omega)$ be a fixed $2n$-dimensional connected symplectic
manifold, with associated Poisson algebra
$\p(M):=\big(C^\infty(M),\{\cdot\, ,\cdot\}\big)$, where
$\{\cdot\, ,\cdot\}$ is the Poisson bracket. We will abbreviate
$\p(M)$ by $\p$ when $M$ is fixed in context.

To start the discussion, we state what it means to ``prequantize'' a
Poisson algebra.

\begin{defn}$\,\,$ {\rm Let $\oo$ be any Poisson
subalgebra\footnote{\,By this we mean a linear subspace of $\p$ which
is closed under Poisson bracket (but not necessarily under
multiplication), i.e., a Lie subalgebra.} of $\p$ containing the
constant function $1$. A {\sl prequantization\/} of $\oo$ is a linear
map
$\q$ from $\oo$ to the linear space Op($D$) of symmetric operators
which preserve a fixed dense domain $D$ in some separable Hilbert
space
$\h$ such that for all $f,g \in \oo$

\begin{enumerate}
\begin{enumerate}
\item[(Q1)] ${\cal Q}\big(\{f,g\}\big) =
\frac{i}{\hbar}\big[{\q}(f),{\q}(g)\big]$,
\vskip 6pt
\item[(Q2)] ${\cal Q}(1) = I$, and
\vskip 6pt
\item[(Q3)] if the Hamiltonian vector field $X_f$ of $f$ is complete,
then
$\q(f)$ is essentially self-adjoint on $D$.
\end{enumerate}
\end{enumerate}

\noi If $\oo = \p$, the prequantization is said to be {\sl full\/}.}
\end{defn}

\noi {\bf Remarks:} {\bf 1.} By virtue of (Q1) a prequantization $\q$
of
$\oo$ is essentially a Lie algebra representation of
$\oo$ by symmetric operators. (More precisely: if we set $\pi(f) =
\frac{i}{\hbar}\,\q(f)$, then $\pi$ is a true Lie algebra
representation by skew-symmetric operators on $D$ equipped with the
commutator bracket. We will blur the distinction between $\pi$ and $\q$.)
In this context there are several additional requirements we could place
upon
$\q$, such as irreducibility and integrability. However, we do not
want to be too selective at this point, so we do not insist on these;
they will be discussed as the occasion warrants.

{\bf 2.} Condition (Q2) reflects the fact that if an observable
$f$ is a constant $c$, then the probability of measuring $f = c$ is
one regardless of which quantum state the system is in. It also
serves to eliminate some ``trivial'' possibilities, such as the
regular representation $f
\mapsto X_f$ on
$L^2(M,\omega^n)$, where $X_f$ is the Hamiltonian vector field of
$f.$

{\bf 3.} Regarding (Q3), we remark that in contradistinction with \vh\
\cite{vH1}, we do not confine our considerations to only those
classical observables whose Hamiltonian vector fields are complete.
Rather than taking the point of view that ``incomplete'' classical
observables cannot be quantized, we simply do not demand that the
corresponding quantum operators be essentially self-adjoint
(``e.s.a.''). We do not imply by this that symmetric operators which
are not e.s.a. are acceptable as physical observables; as is well
known, this is a controversial point.

{\bf 4.}  Notice that no assumptions are made at this stage regarding
the multiplicative structure on $\oo$ vis-\`a-vis $\q$. This is
mainly for historical reasons: in classical mechanics the Lie algebra
structure has played a more dominant role than the ring structure on
$C^{\infty}(M)$, so it is natural to concentrate on the former. This
is also the approach favored by Dirac \cite{di} and the geometric
quantization theorists \cite{So70,Wo}. For more algebraic treatments,
see
\cite{As,Em,vn}. The ring structure is emphasized to a much greater
degree in deformation quantization theory
\cite{Bayen e.a.,ri2,ri3}.

\ms

Prequantizations in this broad sense (even full ones) are usually
easy to construct, cf. \cite{Ch2,Ur,Wo}. Van Hove was the first to
construct a full prequantization of $\p({\bf R}^{2n})$ \cite{vH1}. It
goes as follows: the Hilbert space $\h$ is
$L^2({\bf R}^{2n})$, for $D$ we take the Schwartz space ${\cal S}({\bf
R}^{2n})$ of rapidly decreasing smooth functions (for instance), and
for $f
\in
\p({\bf R}^{2n})$,
\begin{equation}
\q(f) = -i\hbar\sum_{k=1}^n\left[\frac{\partial f}{\partial
p_k}\left(\frac{\partial}{\partial q^k} -\frac{i}{\hbar}\,p_k\right) -
\frac{\partial f}{\partial q^k}\frac{\partial}{\partial p_k}\right] +
f.
\label{eq:preq}
\end{equation}

As luck would have it, however, prequantization representations of the
entire Poisson algebra of a symplectic manifold tend to be flawed.
For example, the Van Hove prequantization of $\p({\bf R}^{2n})$,
when restricted to the Heisenberg subalgebra span$\{1,\,p_i,\,q^i\,|
\,i=1,\ldots,n\}$, is not unitarily equivalent to the Schr\"odinger
representation (which it ought to be, in the context of a particle
moving in
$\r^n$ with no superselection rules)
\cite{bl1,Ch1}. (Recall that the {\sl Schr{\"o}dinger
representation\/} of the Heisenberg algebra
\begin{eqnarray*} {\rm h(2}n{\rm )} & \hspace{-1ex} = \hspace{-1ex} & {\rm
span}\big\{P_i,\,Q^i,\,T,\; i = 1,\ldots,n \,\big|\, [P_j,Q^k] =
\delta^k_j T, \;\; [P_j,P_k] = 0,
\\  & &
\rule{6ex}{0ex} [Q^j,Q^k] = 0,
\;\; [P_j,T] = 0, \;\; [Q^j,T] = 0\big\}
\end{eqnarray*}

\noi is defined to be
\begin{equation} Q^i \mapsto q^i,\;\;\; P_j \mapsto -i\hbar \,
\partial /{\partial q^j},\;\;\; \mbox{and}\;\;\; T
\mapsto I
\label{eq:srep}
\end{equation}

\noi on the domain $\s (\r^n) \subset L^2({\r^{n}}).$ It is
irreducible in the sense given in the next section.) There are
various ways to see this; we give
\vh's original proof
\cite[\S 17]{vH1} as it will be useful later. Take $n=1$ for
simplicity. First, define a unitary operator $F$ on $L^2(\r^2)$ by

\[(F\psi)(p,q) = \frac{1}{\sqrt h}\int^{\infty}_{-\infty}
e^{ipv/\hbar}\psi(v,q-v)\,dv.\]

\noi Then for each fixed $j = 0,1,\ldots$ take
$\h_j$ to be the closure in $L^2(\r^2)$ of the linear span of
elements of the form
$Fh_{jk}$, where $h_{jk}(p,q) = h_j(p)h_k(q)$ and $h_j$ is the
Hermite function of degree $j$. Now from (\ref{eq:preq}),

\[\q(q) = i\hbar \, \frac{\partial}{\partial p} + q,\;\;\;\q(p) =
-i\hbar
\,\frac{\partial}{\partial q}.\]

\noi These operators are e.s.a.\   on
$\s (\r^2),$ and one may verify that they strongly com\-mute with the
orthogonal projectors onto the closed subspaces
$\h_j$.\footnote{\,Recall that two e.s.a.\ (or, more generally,
normal) operators {\sl strongly commute} iff their spectral
resolutions commute, cf.
\cite[\S VIII.5]{ReSi}. Two operators $A,\,B$ {\sl weakly commute} on
a domain $D$ if they commute in the ordinary sense, i.e., $[A,B]$ is
defined on $D$ and vanishes.} Thus the
\vh\ prequantization of $\p(\r^2)$ is reducible when restricted to the
Heisenberg subalgebra and hence does not produce the Schr\"odinger
representation. Moreover the association $Fh_{jk}(p,q) \mapsto
c_{j}h_k(q),$ where the $c_{j}$ are normalization constants, provides
a unitary equivalence of each subrepresentation of h(2) on
$\h_j$ with the Schr\"odinger representation on $L^2(\r)$,  from
which we see that the multiplicity of the latter is infinite in the
\vh\ representation.

Likewise, the Kostant-Souriau prequantizations of
$S^2$ do not reproduce the familiar spin representations of the
unitary algebra u(2). We realize $S^2$ as a coadjoint orbit of SU(2)
in $\r^3$ according to ${\bf S\cdot S} = s^2$, where ${\bf S} =
(S_1,S_2,S_3)$ is the spin vector and
$s > 0$ is the classical spin. It comes equipped with the symplectic
form
\begin{equation}
\omega=\frac{1}{2s^2}\sum_{i,j,k = 1}^3\epsilon_{ijk}\,S_i\,
dS_j\wedge dS_k.
\label{eq:sfs2}
\end{equation}

\noi Now the de Rham class $[\,\omega/h\,]$ is integral iff $s =
\frac{n}{2}\hbar$, where
$n$ is a positive integer, and the corresponding Kostant-Souriau
prequantization line bundles can be shown to be
$L^{\otimes n}$ where $L$ is the dual of the universal line bundle
over
$S^2$ \cite{ki}. The corresponding prequantum Hilbert spaces $\h_n$
can thus be identified with spaces of square integrable sections
$\psi$ of these bundles w.r.t. the inner product
\[\langle \psi,\phi \rangle = \frac{i}{2\pi}
\int_{\scriptscriptstyle {\bf C}}\frac{\overline{\psi(z)}\phi(z)\,dz
\wedge d{\bar z}}{(1+z{\bar z})^{n+2}}\]

\noi where $z = (S_1 + iS_2)/(s - S_3)$, cf. \cite{Wo}. But these
$\h_n$ are infinite-dimensional, whereas the standard representation
spaces for quantum spin $s = \frac{n}{2}\hbar$ have dimension $n+1$.

In both examples the prequantization Hilbert spaces are ``too big.''
The main problem is how to remedy this, in other words, how to modify
the notion of a prequantization so as to yield a genuine
${\it quantization.}$

It is here that the ideas start to diverge, because there is less
agreement in the literature as to what constitutes a quantization
map. Some versions define it as a prequantization, not necessarily
defined on the whole of
$\p$, which is irreducible on a ``basic set''
$\bb\subset\p$ \cite{ki}. This is in line with the group theoretical
approach to quantization \cite{Is}, in which context $\bb$ is
identified with the Lie algebra of a symmetry group;\footnote{\,We
typically identify an abstract Lie algebra with its isomorph in
$\p$.} quantization should then  yield an irreducible representation
of this algebra.  For example, when
$M=\r^{2n}$, one usually takes $\bb$ to be the Heisenberg algebra
${\rm h(2{\em n})} = {\sp}\{1,\, p_i,\, q^i\,|\,i=1,\ldots, n\}$ of
polynomials of degree at most one. Similarly, when
$M=S^2$, one takes for $\bb$ the unitary subalgebra ${\rm
u(2)}={\sp}\{1,\, S_1,\, S_2,\, S_3\}$ of spherical harmonics of
degree at most one, where
$S_i$ are the spin generators. We will plumb in detail the rationale
behind these choices of $\bb$ in the next section.

A different approach to quantization is to require a prequantization
$\q$ to satisfy some ``Von~Neumann rule,'' that is, some given
relation between the classical multiplicative structure of $\p$ and
operator multiplication on
$\h$. (Note that thus far in our discussion the multiplication on
$\p$ has been ignored, and it is reasonable to require that
quantization preserve at least some of the ring structure of $\p$,
given that the Leibniz rule intertwines pointwise multiplication with
the Poisson bracket.) There are many different types of such rules
\cite{Co,Fo,KLZ,KS,Ku,MC,vn}, the simplest being of the form:
\begin{equation}
\q(\varphi\circ f)=\varphi\big(\q(f)\big)
\label{eq:vnr}
\end{equation}

\noi for some distinguished observables
$f\in\p$, and certain smooth functions
$\varphi\in C^\infty(\r)$. (Technically, if $\varphi$ is not a
polynomial, then
$\q(f)$ must be e.s.a.\  for $\varphi\big(\q(f)\big)$ to be defined.)
In the case of
$M=\r^{2n}$, \vn\ states that the physical interpretation of the
quantum theory requires (\ref{eq:vnr}) to hold for {\em all\/} $f \in
\p$ and $\varphi\in C^\infty(\r)$
\cite{vn}. However, it is easy to see that this is impossible (simple
demonstrations are given in
\cite{a-b,Fo} as well as \S 4.1 following); thus the qualifiers in the
definition above. In this example, one typically ends up imposing the
squaring
\vn\ rule
$\varphi(x) = x^2$ on elements $x$ of h(2$n$). The relevant rules for
the sphere turn out to be somewhat less intuitive; they take the form
$\q\big(S_i\,\!^2\big) = a\q(S_i)^2 + cI$ for $i = 1,\,2,\,3$, where
$a,c$ are undetermined (representation-dependent) constants subject
only to the constraint that $ac \neq 0$. Derivations of these rules
in these two examples are given in \S4 and \cite{GGH}.

Another type of quantization is obtained by ``polarizing'' a
prequantization representation
\cite{Wo}. Following Blattner \cite{bl1}, we paraphrase it
algebraically as follows. Start with a {\sl polarization\/}, i.e., a
maximally commuting Poisson subalgebra $\A$ of $\p$. Then require for
the quantization map
$\q$ that the image $\q(\A)$ be ``maximally commuting'' as operators.
If
$\q(\A)$ consists of bounded operators, this means that the weak
operator closure of the *-algebra generated by $\q(\A)$ \big($=
\q(\A)''$\big) is maximally commuting in $B(\h )$. If
$\q(\A)$ contains unbounded operators, one should look for a
generating set of normal operators in $\q(\A)$, and require that the
Von~Neumann algebra generated by their spectral projections is
maximally commuting. One can then realize the Hilbert space
$\h $ as an $L^2$-space over the spectrum of this \vn\ algebra on
which this algebra acts as multiplication operators. There will also
be a cyclic and separating vector for such an algebra, which provides
a suitable candidate for a vacuum vector. Thus another motivation for
polarizations is that a maximally commuting set of observables
provides a set of compatible measurements, which can determine the
state of a system. When $M=\r^{2n}$, one often takes the ``vertical''
polarization $\A =
\big\{f(q^1,\ldots,q^n)\big\}$, in which case one
recovers the usual position or coordinate representation. However, in
some instances, such as $S^2$, it is useful to broaden the notion of
polarization to that of a maximally commuting subalgebra of the {\em
complexified\/} Poisson algebra
$\p_{\scriptscriptstyle {\bf C}}.$ Then, thinking of $S^2$ as $\C
P^1$, we may take the ``antiholomorphic'' polarization $\A =
\{f(z)\},$ which leads to the usual representations for spin.
For treatments of polarizations in the context of deformation
quantization, see
\cite{Fr,He}.

Thus, informally, a ``quantization'' could be defined as a
prequantization which incorporates one (or more) of the three
additional requirements above (or possibly even others). Before
proceeding, however, there are two points we would like to make.

The first is that it is of course not enough to simply state the
requirements that a quantization map should satisfy; one must also
devise methods for implementing them in examples. Thus geometric
quantization theory, for instance, provides a specific technique for
polarizing certain (Kostant-Souriau) prequantization representations
\cite{bl1,ki,So70,Wo}. However, as we are interested here in the
structural aspects of quantization theory, and not in specific
quantization schemes, we do not attempt to find such implementations.

Second, these three approaches to a quantization map are not
independent; there exist subtle connections between them which are
not well understood. For instance, demanding that a prequantization
be irreducible on some basic set typically leads to the appearance of
\vn\ rules; this is how the \vn\ rules for $\r^{2n}$ and $S^2$
mentioned above arise. We will delineate these connections in
specific cases in
\S 4, and more generally in \S 5.

\ms

At the core of each of the approaches above is the imposition -- in
some guise -- of an irreducibility requirement, which is used to ``cut
down'' a prequantization representation. Since this is most apparent
in the first approach, we will henceforth concentrate on it. We will
tie in the two remaining approaches as we go along.

So let $\oo$ be a Poisson subalgebra of $\p$, and suppose that $\bb
\subset \oo$ is a ``basic set'' of observables. Provisionally, we
take a {\sl quantization} of the pair $({\oo},{\bb})$ to mean a
prequantization $\q$ of $\oo$ which (among other things) irreducibly
represents $\bb$. In the next section we will make this more precise,
as well as examine in detail the criteria that $\bb$ should satisfy.

\ms

Natural issues to address for quantizations are existence, uniqueness
and classification, and functoriality. For {\em pre\/}quantizations
these questions already have partial answers from geometric
quantization theory. So for instance we know that if $(M,\omega)$
satisfies the integrality condition
$[\,\omega/h\,]\in H^2(M,\z)$, then full prequantizations of the
Poisson algebra $\p$ exist, and that certain types of these -- the
Kostant-Souriau prequantizations -- can be  classified
cohomologically  \cite{Ur,Wo}. For some limited types of manifolds
the functorial properties of these prequantizations were considered
by Blattner
\cite{bl1}. However, as there are prequantizations not of the
Kostant-Souriau type \cite{av,Ch2}, these questions are still open in
general (especially for manifolds which violate the integrality
condition
\cite{We}).

For quantization maps these questions are far more problematic. Our
main focus will be on the existence of {\sl full quantizations}, by which
we mean a quantization of $(\p,\bb)$ for some appropriately chosen
basic set $\bb$. As indicated earlier, this is poorly understood. In
our terminology, the classical (strong) result of Groenewold states
that there is no full quantization of
$\big(\p(\r^{2n}),{\rm h}(2n)\big)$, while the more recent result of
\cite{GGH} implies essentially the same for $\big(\p(S^2),{\rm
u}(2)\big)$. On the other hand, nontrivial full quantizations do
exist: one can construct such a  quantization of $T^2$ \cite{Go}.
However, it does seem that nonexistence results are the rule. In the
absence of a full quantization, then, it is important to determine
the maximal subalgebras
$\oo$ of $\p$ for which $(\oo,\bb)$ can be quantized. This we will
investigate for
$\r^{2n}$ and $S^2$ in \S\S 4 and 5. At present, questions of
uniqueness and classification can only be answered in specific
examples.

\end{section}


\begin{section}{Basic Sets and Quantization}

Our first goal here is to make clear what we mean by a basic set of
observables
$\bb\subset\p$. Such sets, in one form or another, play an important
role in many quantization methods, such as geometric quantization
\cite{ki}, deformation quantization \cite{Bayen e.a.,Fr} and also the
group theoretic approach \cite{Is}.

We start with the most straightforward case, that of an ``elementary
system'' in the terminology of Souriau
\cite{So70,Wo}. This means that $M$ is a homogeneous space for a
Hamiltonian action of a finite-dimensional Lie group $G$. The appeal
of an elementary system is that it is a classical version of an
irreducible representation: using the transitive action of
$G$, one can obtain any classical state from any other one, in direct
analogy with the fact that every nonzero vector in a Hilbert space
$\h$ is cyclic for an irreducible unitary representation of $G$ on
$\h$ \cite[\S 5.4]{b-r}. Now notice that the span
$\cal J$ of the components of the associated (equivariant) momentum
map satisfies:
\begin{enumerate}
\begin{enumerate}
\item[(J1)] $\jj$ is a finite-dimensional Poisson subalgebra of $\p$,
\vskip 6pt
\item[(J2)] the Hamiltonian vector fields
$X_f,{f\in\jj}$, are complete, and
\vskip 6pt
\item[(J3)]  $\{X_f \,|\,{f\in\jj}\}$ span the tangent spaces to $M$
everywhere.
\end{enumerate}
\end{enumerate}

\noi  For both $M=\r^{2n}$ and $S^2$, the basic sets are precisely of
this type: from the elementary systems of the Heisenberg group
H(2$n$) acting on
$\r^{2n}$, and the unitary group U(2) acting on $S^2$, we have for
$\jj$ the spaces span$\{1,\,p_i,\,q^i\,|\,i = 1,\ldots,n\}$ and
span$\{1,\,S_1,\,S_2,\,S_3\}$, respectively.

Property (J3) is just an infinitesimal restatement of transitivity,
and so we call a subset $\bb$ of $C^{\infty}(M)$ {\sl transitive\/}
if it satisfies this condition. Kirillov \cite{ki} uses the
terminology ``complete,'' motivated by the fact that such a set of
observables locally separates classical states. In this regard, the
finite-dimensionality criterion in (J1) plays an important role
operationally: it guarantees that a {\em finite\/} number of
measurements using this collection of observables will suffice to
distinguish any two nearby states.

A subset $\bb \subset \p$  satisfying (J1)--(J3) is a prototypic basic
set. However, there need not exist basic sets in this sense for
arbitrary $M$. For instance, if
$M = T^2$, the self-action of the torus is not Hamiltonian -- there
is no momentum map -- and consequently it is difficult to isolate
such a basic set. Thinking of $T^2$ as
$\r^{2}/\z^{2}$, a natural choice for $\bb$ would be

\[{\cal B} = \mbox{\rm span}\{1,\sin 2\pi x,\cos 2\pi x,\sin 2\pi
y,\cos 2\pi y\},\]

\noi but this is not a subalgebra. One alternative would be to
consider instead the Poisson algebra $\wp(\bb)$ generated by $\bb$.
However, this algebra (viz. the set of trigonometric polynomials) is
infinite-dimensional, and there is nothing gained operationally in
using an infinite-dimensional algebra to distinguish classical states.
Furthermore, as will be shown below, other problems arise if one
insists that $\bb$ always be a subalgebra of $\p.$ We will therefore
retain the finite-dimensionality assumption, but merely require that
$\bb$ be a linear sub{\em space\/} as opposed to a Poisson sub{\em
algebra} of
$\p$. Thus we make:

\begin{defn} {\rm \, A {\sl basic set of observables\/} $\bb$ is a
linear subspace of $\p(M)$ such that:
\vspace{-6pt}
\begin{enumerate} \begin{enumerate}
\item[(B1)] $\bb$ is finite-dimensional,
\vskip 6pt
\item[(B2)] the Hamiltonian vector fields
$X_f,{f\in\bb}$, are complete,
\vskip 6pt
\item[(B3)] $\bb$ is transitive,
\vskip 6pt
\item[(B4)] $1 \in \bb$, and
\vskip 6pt
\item [(B5)] $\bb$ is a minimal space satisfying these requirements.
\end{enumerate} \end{enumerate}}
\end{defn}

We spend some time elaborating on this definition. First,
condition (B4) is algebraically natural, as discussed in Remark {\bf
2}. (This also explains why, for $S^2$, we take $\bb = {\rm u}(2)$
rather than su(2).) Second, the minimality condition (B5) is crucial.
{}From a physical or operational point of view, it is not obvious
that this is necessary, as long as $\bb$ is finite-dimensional. But
the quantization of a pair
$(\oo,\bb)$ with $\bb$ non-minimal can lead to physically incorrect
results.

Here is an example of this phenomenon. First observe that the
extended symplectic group HSp(2$n$,$\r$) (which is the semi-direct
product of the symplectic group Sp(2$n$,$\r$) with the Heisenberg
group H(2$n$)) acts transitively on
$\r^{2n}$. This action has a momentum map whose components consist of
all inhomogeneous quadratic polynomials in the $q^i$ and
$p_i$. The corresponding subalgebra $\jj$ satisfies all the
requirements for a basic set save minimality. Now consider again the
Van Hove prequantization $\q$ of
$\p(\r^{2n})$ for $n=1$. In
\cite[\S 17]{vH1} it is shown that $\q$ is completely reducible when
restricted to the hsp(2$,\r$)-subalgebra
$\jj$. In fact, there exist exactly two nontrivial
HSp(2,$\r$)-invariant closed subspaces $\h_{\pm}$ in $L^2(\r^{2})$,
namely (cf. \S 2)

\[\h_+ = \bigoplus_{j \;{\rm even}}\h_j\;\;\;{\rm and}\;\;\; \h_- =
\bigoplus_{j
\;{\rm odd}}\h_j.\]

\noi If we denote the corresponding subrepresentations of
$\jj$ on
${\cal S}(\r^{2}) \cap \h_{\pm}$ by
$\q_{\pm}$, then it follows that $\q_{\pm}$ are quantizations of the
pair
$(\jj,\jj).$ But these quantizations are physically unacceptable,
since -- just like the full prequantization $\q$ -- they are
reducible when further restricted to h(2)
$\subset$ hsp(2$,\r$). On the one hand, asking for a quantization of
$(\jj,\jj)$ in this context is clearly the wrong thing to do, since
compatibility with Schr\"odinger quantization devolves upon the
irreducibility of an h(2$n$) subalgebra, not an hsp(2$n,\r$) one. But
on the other, this example does illustrate our point.

As well, the minimality requirement reinforces our assertion that it
will not do to demand that the basic set be a Poisson subalgebra
rather than a linear subspace. For consider again the torus and
define the basic sets
\[{\cal B}_k = \mbox{\rm span}\{1,\sin 2\pi kx,\cos 2\pi kx,\sin 2\pi
ky,\cos 2\pi ky\}\]

\noi for $k = 1,2,\ldots$ Each $\bb_k$ is a minimal transitive
subspace. But
$\wp(\bb_k)$ is {\em not\/} a minimal transitive subalgebra for any
$k$, since $\wp(\bb_k) \supset \wp(\bb_{2k}) \supset
\cdots\,$. In fact, there probably does not exist a minimal
transitive Poisson sub{\em algebra\/} of
$\p(T^2)$.


Finally, we consider the transitivity requirement (B3). While (B3) is
geometrically natural, there are other conditions one might use in
place of it. By way of motivation, consider a unitary representation
$U$ of a Lie group
$G$ on a Hilbert space $\h$. The representation
$U$ is irreducible iff the *-algebra ${\cal U}$ of bounded operators
generated by $\{U(g)\,|
\,g \in G\}$ is irreducible, in which case we have the following
equivalent characterizations of irreducibility:

\begin{enumerate}\begin{enumerate}
\item[(I1)] The commutant ${\cal U}'= \C I,$ and
\vskip 6pt
\item[(I2)] the weak operator closure of $\cal U$ is the algebra of
all bounded operators: $\overline{\cal U}^w=B(\h)\quad(={\cal U}'').$

\end{enumerate}\end{enumerate}

\noi That (I1) is equivalent to irreducibility is the content of
Schur's Lemma. Property (I2) means that all bounded operators can be
built from those in
$\cal U$ by weak operator limits. It follows from (I1), the \vn\
density theorem \cite[Cor. 2.4.15]{BR}, and the fact that
${\cal U}'=\big(\,\overline{\cal U}^w\big)'$. Clearly (I2) implies
(I1).

These restatements of irreducibility have the following classical
analogues for a set ${\cal F}$ of observables:

\begin{enumerate}\begin{enumerate}
\item[(C1)] $\{f,g\}=0$ for all
$f\in {\cal F}$ implies
$g$ is constant, and
\vskip 6pt
\item[(C2)] ${\cal F}$ generates a dense subspace in
$C^{\infty}(M)$ under linear combinations and pointwise
multiplication.

\end{enumerate}\end{enumerate}

\noi For (C2) a topology on $C^{\infty}(M)$ must be decided on, and
we will use the topology of uniform convergence on compacta of a
function as well as its derivatives.

Because the algebraic structures of classical and quantum mechanics
are different,  (C1) and (C2) lead to inequivalent notions of
``classical irreducibility.'' It is not difficult to verify that (C1)
$\Leftarrow$ (B3) $\Leftarrow$ (C2) strictly. In principle either of
(C1) or (C2) could serve in place of (B3). Indeed, since on
$C^{\infty}(M)$ one has two algebraic operations, it is natural to
consider irreducibility in either context: in terms of the
multiplicative structure (C2), or the Poisson bracket (C1). However,
it turns out that (C1) is too weak for our purposes, while (C2) is too
strong.

The nondegeneracy condition (C1) is equivalent to the statement that
observables in $\bb$ locally separate states {\em almost\/}
everywhere \cite{ki}. It is also implied by the statement that the
Hamiltonian vector fields of elements of $\bb$ span the tangent
spaces to
$M$ almost everywhere. Consequently, it would not do to replace (B3)
by (C1) in the definition of basic set because, e.g., in the case of
the sphere span$\{1,\,S_1,\,S_2,\,S_3\}$ would no longer be minimal,
which seems both awkward and unreasonable. The same is true for the
sets
$\bb_k$ on the torus, as well as the algebra sp(2$n,\r) \times \r$ on
the homogeneous space $\r^{2n}\setminus \{{\bf 0}\}$. Condition (C2)
is satisfied for the unitary algebra on
$S^2$, the Heisenberg algebra on
$\r^{2n}$ and the set $\bb_1$ on $T^2$. But it fails for the
symplectic algebra on $\r^{2n}\setminus \{{\bf 0}\}$ -- since the
subspace generated by sp$(2n,\r)$ consists of even functions, and
also for the sets
$\bb_k$ with $k>1$ on the torus -- since the subspace they generate
consists of doubly periodic functions of period $k>1$. On the other
hand, all these examples satisfy the transitivity requirement (B3),
which shows that this is a reasonable condition to impose.

There is still no guarantee in general that sets satisfying (B1)-(B5)
will exist for a given symplectic manifold $M$. However, we can
satisfy all conditions except (B2) as follows. Choose an embedding $M
\ra
\r^K$ for some sufficiently large $K$, and let $\{f^1,\ldots,f^K\}$
be the restrictions of the standard coordinates on $\r^K$ to $M$.
Then ${\cal F} =$ span$\{1,f^1,\ldots,f^K\}$ satisfies (B1), (B3) and
(B4). If ${\cal F}$ is not minimal as it stands, one may discard
elements of this set until what remains is minimal. If we allow
infinite-dimensional basic sets, then  localizing ${\cal F}$ by means
of a partition of unity enables us to satisfy (B2), but now of course
(B1) has to be abandoned. In any event, we do know that basic sets
will exist whenever ({\em i\/}) $M$ is a Hamiltonian homogeneous
space, or ({\em ii\/}) $M$ is compact.

Other properties that basic sets might be required to satisfy are
discussed in \cite{Is}. For our purposes, (B1)-(B5) will suffice.

\ms

We are now ready to discuss what we mean by a ``quantization.'' Let
$\oo$ be a Poisson subalgebra of $\p$, and suppose that $\bb \subset
\oo$ is a basic set of observables. Two eminently reasonable
requirements to place upon a quantization are irreducibility and
integrability \cite{b-r,fl,Is,ki}.

Irreducibility is of course one of the pillars of the quantum theory,
and we have already seen the necessity of requiring that quantization
represent
$\bb$ irreducibly. We must however be careful to give a precise
definition since the operators
$\q(f)$ are in general unbounded (although, according to (B2) and
(Q3), all elements of
$\q(\bb)$ are e.s.a.). So let
${\cal X}$ be a set of e.s.a.\  operators defined on a common
invariant dense domain $D$ in a Hilbert space $\h$. Then
${\cal X}$ is {\sl irreducible\/} provided the only bounded
self-adjoint operators which strongly commute with all $X \in {\cal
X}$ are multiples of the identity. While this definition is fairly
standard, and well suited to our needs, we note that other notions
of irreducibility can be found in the literature \cite{b-r,mmsv}.

Given such a set ${\cal X}$ of operators, let $\uu ({\cal X})$ be the
$^*$-algebra generated by the unitary operators $\big\{\exp(it \ol X)
\,|\, t \in \r,\; X \in {\cal X}$\big\}, where $\ol X$ is the closure
of
$X$. Then by Schur's Lemma ${\cal X}$ is irreducible iff the only
closed subspaces of
$\h$ which are invariant under $\uu ({\cal X})$ are $\{0\}$ and $\h$.

Turning now to integrability, we first consider the case when the
Poisson algebra $\wp(\bb)$ generated by the basic set is
finite-dimensional. Then it is natural to demand that the Lie algebra
representation
$\q\big(\wp(\bb)\big)$ on $D
\subset \h$ be {\sl integrable\/} in the following sense:  there
exists a unitary representation ${\mit \Pi}$ of some Lie group with
Lie algebra (isomorphic to) $\wp(\bb)$ on $\h$ such that
$d{\mit \Pi}(f)|D = \q(f)$ for all $f \in \wp(\bb)$, where $d{\mit
\Pi}$ is the derived representation of ${\mit \Pi}.$ For this it is
neither necessary {\em nor\/} sufficient that elements of
$\wp(\bb)$ quantize to e.s.a. operators on $D$.\footnote{\,See Remark
{\bf 8} following.} But integrability will follow from the following
result of Flato et al., cf. \cite{fl} and \cite[Ch. 11]{b-r}.
\begin{prop} Let ${\cal G}$ be a real finite-dimensional Lie algebra,
and let $\pi$ be a representation of ${\cal G}$ by skew-symmetric
operators on a common dense invariant domain $D$ in a Hilbert space
$\h$. Suppose that $\{\xi_1,\ldots,\xi_k\}$ generates ${\cal G}$ by
linear combinations and repeated brackets. If
$D$ contains a dense set of separately analytic vectors for
$\big\{\pi(\xi_1),\ldots,\pi(\xi_k)\big\}$, then there exists a unique
unitary representation ${\mit \Pi}$ of the connected simply connected
Lie group with Lie algebra ${\cal G}$ on $\h$ such that
$d{\mit \Pi}(\xi)|D = \pi(\xi)$ for all $\xi \in {\cal G}$.
\label{prop:fs}
\end{prop}

(We recall that if $X$ is an operator on $\h$, a vector $\psi$ is
{\sl analytic\/} for $X$ provided the series
\[\sum_{k=0}^{\infty}\frac{\|X^k\psi\|}{k!}\,t^k\]

\noi is defined and converges for some $t>0.$ If
$\{X_1,\ldots,X_k\}$ is a  set of operators
defined on a common invariant dense domain $D$, a vector $\psi
\in D$ is {\sl separately analytic\/} for $\{X_1,\ldots,X_k\}$ if
$\psi$ is analytic for each $X_j$. By a slight abuse of terminology,
we will say that a vector is separately analytic for a linear
space of operators ${\cal X}$ if it is separately analytic for some
basis $\{X_1,\ldots,X_k\}$ of ${\cal X}$.)

However, it may happen that $\wp(\bb)$ is not finite-dimensional (as
in the case of the torus). Then there need not exist an
(infinite-dimensional) Lie group having $\wp(\bb)$ as its Lie
algebra. Even if such a Lie group existed, integrability is far from
automatic, and technical difficulties abound. Thus we will not insist
that a quantization be integrable in general. On the other hand, the
analyticity requirement in Proposition~\ref{prop:fs} makes sense
under all circumstances,\footnote{\,As long as $\wp(\bb)$ is finitely
generated, which is assured by (B1).} and does guarantee
integrability when
$\wp(\bb)$ is ``nice,'' so we will adopt it in lieu of integrability.

Therefore we have at last:

\begin{defn} \, {\rm  A {\sl quantization\/} of the pair
$({\oo},{\bb})$ is a prequantization $\q$ of $\oo$ on Op($D$)
satisfying
\begin{enumerate}
\begin{enumerate}
\item[(Q4)] $\q(\bb) = \left\{\q(f)\,|\,f \in \bb\right\}$ is an
irreducible set, and
\vskip 6pt
\item[(Q5)] $D$ contains a dense  set of separately analytic vectors
for
$\q(\bb).$
\label{def:q}
\end{enumerate}
\end{enumerate}}

\noi {\rm A quantization $\q$ is {\sl nontrivial\/} if the
representation space is neither one- nor zero-dimensional.}
\label{defn:quant}
\end{defn}

\noi {\bf Remarks: 5.}  There are a number of analyticity assumptions
similar to (Q5) that one could make \cite{fl}; we have chosen the
weakest possible one.

{\bf 6.} (Q5) is not a severe restriction: when $\wp(\bb)$ is
finite-dimensional, it is always possible to find representations
of it on domains $D$ which satisfy this property \cite{fl}. On the
other hand, nonintegrable representations do exist in general
\cite[p. 247]{fl}.

{\bf 7.} Proposition~\ref{prop:fs} requires that a specific
generating set for $\q\big(\wp(\bb)\big)$ be singled out. This also
is not a severe restriction: in examples, $\bb$ is usually given as
the linear span of such a set. It is possible that (Q5) could be
satisfied for one such set but not another, but Remark {\bf 6} shows
that the domain
$D$ can be chosen in such a way that this cannot happen if $\wp(\bb)$
is finite-dimensional.

{\bf 8.} With regard to essential self-adjointness vis-\`a-vis
integrability, we observe that the operators $\q(f)$ for $f
\in \wp(\bb)$ need {\em not\/} be e.s.a.\ on $D$;
Proposition~\ref{prop:fs} only guarantees that they have s.a.\
extensions. This is consistent with the fact that the Hamiltonian
vector fields of elements of $\wp(\bb)$ are not necessarily complete.
When $f$ is complete (e.g., $f \in
\bb$), (Q3) requires that $\q(f)$ be e.s.a.\ and the proposition
correspondingly yields $d{\mit \Pi}(f) = \ol {\q(f)}.$ That essential
self-adjointness alone is not sufficient to guarantee integrability
is well known, cf. \cite[\S VIII.5]{ReSi}.

{\bf 9.} It is important to note that irreducibility does not imply
integrability. For instance, there is an irreducible representation
of h(2) which is not integrable \cite[p. 275]{ReSi}.

\ms

We end this section with a brief comment on the domains $D$ appearing
in Definition~\ref{def:q}. For a representation $\pi$ of a Lie
algebra ${\cal G}$ on a Hilbert space $\h$, there is typically a
multitude of common, invariant dense domains that one can use as
carriers of the representation. (See \cite[\S 11.2]{b-r} for a
discussion of some of the possibilities.) But what is ultimately
important for our purposes are the closures $\ol{\pi(\xi)}$ for $\xi
\in {\cal G}$, and not the
$\pi(\xi)$ themselves. So we do not want to distinguish between two
representations $\pi$ on Op$(D)$ and $\pi'$ on Op$(D')$ whenever
$\ol{\pi(\xi)} = \ol{\pi'(\xi)}$, in which case we say that $\pi$ and
$\pi'$ are {\sl coextensive}. In particular, it may happen that the
given domain $D$ for a representation $\pi$ does not satisfy (Q5),
but there is an extension to a coextensive representation $\pi'$ on a
domain $D'$ that does.\footnote{\, A simple illustration is provided
by the Schr\"odinger representation (\ref{eq:srep}) with $D =
C^{\infty}_0(\r^n)$ and $D' =\s (\r^n)$.} In such cases we will
suppose that the representation has been so extended.

\end{section}


\begin{section}{Examples}

In this section we present the gist of the arguments that there are no
nontrivial quantizations of either $\left(\p(\r^{2n}) ,{\rm h}(2n)
\right)$ or $\left(\p(S^2),{\rm u}(2)\right)$. The complete proofs
can be found in \cite{a-m,Ch1,Fo,go80,Gr,GS,vH1,vH2} for $\r^{2n}$ and in
\cite{GGH} for $S^2$. In both cases the  detailed structure of the
Poisson algebras and their representation theory is used, which makes
it hard to generalize these results to other symplectic manifolds.
Finally, we show following \cite{Go} that there is a full
quantization of
$\left(\p(T^2),\bb_k\right)$ for each integer $k>0$, where $\bb_k$ is
the basic set defined in the last section. We also take this
opportunity to point out a defect in the standard presentations of the
\gr-\vh\ theorem for
$\r^{2n}.$


\begin{subsection}{$\r^{2n}$}

Before proceeding with the no-go theorem for
$\r^{2n}$, we remark that already at a purely mathematical level one
can observe a suggestive structural mismatch between the classical
and the quantum formalisms. Since a prequantization is essentially a
Lie algebra representation, it ``compares'' the Poisson algebra
structure of
$\p(\r^{2n})$ with the Lie algebra of  (skew-) symmetric operators
(preserving a dense domain
$D$) equipped with the commutator bracket. But if we take $P
\subset \p(\r^{2n})$ to be the subalgebra of polynomials, Joseph
\cite{Jo} has shown that $P$ has outer derivations, but the
enveloping algebra of the Heisenberg algebra h(2$n$) -- and hence
that of the Schr\"odinger representation thereof on $L^2(\r^n)$ --
has none.


Furthermore, one can see at the outset that it is impossible for a
prequantization to satisfy the ``product $\ra$ anti-commutator'' rule.
Taking $n=1$ for simplicity, suppose $\q$ were a prequantization of
the polynomial subalgebra $P \subset \p(\r^{2})$ for which
\begin{equation}  \q(fg) = {\textstyle \frac{1}{2}}\big(\q(f)\q(g) +
\q(g)\q(f)\big)
\label{eq:ac}
\end{equation}

\noi for all $f,g \in P.$ Take $f(p,q) = p$ and $g(p,q) = q$. Then
\begin{eqnarray*} {\textstyle \frac{1}{4}}\big(\q(p)\q(q) +
\q(q)\q(p)\big)^2 \hspace{-1ex} & = & \hspace{-1ex}
\q(pq)^2 \\ & = & \hspace {-1ex} \q\big(p^2q^2\big)\; = \; {\textstyle
\frac{1}{2}}\big(\q(p)^2\q(q)^2 +
\q(q)^2\q(p)^2\big).
\end{eqnarray*}

\noi Now by (Q1) we have $[\q(p),\q(q)] = -i\hbar$, so  that the
L.H.S. reduces to
\[\q(q)^2\q(p)^2 - 2i\hbar\q(q)\q(p) - {\textstyle \frac{1}{4}}\hbar^2 I\]

\noi while the R.H.S. becomes
\[\q(q)^2\q(p)^2 - 2i\hbar\q(q)\q(p) - \hbar^2 I.\]

\noi As the product $\ra$ anti-commutator rule is equivalent to the
squaring \vn\ rule $\q\big(f^2\big) =\q(f)^2$, this contradiction
also shows that the latter is inconsistent with prequantization. Note
that the contradiction is obtained on quartic polynomials; there is
no problem if consideration is limited to observables which are at
most cubic.

This argument only used axiom (Q1) in the specific instance
$[\q(p),\q(q)] = -i\hbar I$. Consequently, one still obtains a
contradiction if one drops (Q1) and instead insists that $\q$ be
consistent with Schr\"odinger quantization (in which context this one
commutation relation remains valid, cf.~(\ref{eq:srep})). This {\em
manifest\/} impossibility of satisfying the product $\ra$
anti-commutator rule while being consistent with Schr\"odinger
quantization is one reason we have decided to concentrate on the
Lie structure as opposed to the multiplicative structure of
$C^{\infty}(M)$. See \cite{a-b} for further results in this direction.

\ms

We now turn to the no-go theorem for $\r^{2n}$. We shall state the
main results for $\r^{2n}$ but, for convenience, usually prove them
only for
$n=1$. The proofs for higher dimensions are immediate generalizations
of these. In what follows $P$ denotes the subalgebra of polynomials,
$P^k$ the subspace of polynomials of degree at most $k$ and $P_k$ the
space of homogeneous polynomials of degree $k$. Note that
$P^1
\cong {\rm h}(2n)$, $P_2 \cong {\rm sp}(2n,\r)$, and $P^2 \cong {\rm
hsp}(2n,\r)$, where the latter is the Lie algebra of the extended
symplectic group, cf. \S 3.

There are actually several versions of the \gr-\vh\ no-go theorem,
depending upon the properties one wants a quantization to satisfy. We
begin with the weakest result, which requires no assumptions on $\q$
beyond those given in Definition~\ref{defn:quant}.

We first observe that there {\em does\/} exist a quantization
$d\varpi$ of the pair
$(P^2,P^1)$. For
$n=1$ it is given by the familiar formul{\ae}
\begin{eqnarray} d\varpi(q)  = q, & d\varpi(1) = I, & d\varpi(p) =
-i\hbar
\frac{\partial}{\partial q}
\label{eq:p1} \\ d\varpi(q^2) = q^2, & d\varpi(pq) = -i\hbar
{\displaystyle
\left (q \frac{\partial}{\partial q} + \frac{1}{2}\right )}, &
d\varpi(p^2) = -\hbar^2 \frac{\partial^2}{\partial q^{2}}
\label{eq:p2}
\end{eqnarray}

\noi on the domain $\s (\r) \subset L^2(\r).$ Properties (Q1)--(Q3)
are readily verified. (Q4) follows automatically since the
restriction of
$d\varpi$ to
$P^1$ is just the Schr\"od\-inger representation. For (Q5), we recall
the well-known fact that the Hermite functions form a dense set of
separately analytic vectors for $d\varpi(P^1)$. Since these functions
are also separately analytic vectors for $d\varpi(P_2)$
\cite[Prop.~4.49]{Fo}, the operator algebra
$d\varpi(P^2)$ is integrable to a representation $\varpi$ of the
universal cover ${\widetilde{\rm{HSp}}}(2n,\r)$ of
HSp$(2n,\r)$\footnote{\,This representation actually drops to the
double cover of HSp$(2n,\r)$, but we do not need this fact here.}
(thereby justifying our notation ``$d\varpi$'').
$\varpi$ is known as the ``extended metaplectic representation'';
detailed discussions of it may be found in \cite{Fo,GS}.

We call $d\varpi$ the ``extended metaplectic quantization.'' It has
the following crucial property.
\begin{prop} The extended metaplectic quantization is the {\em
unique\/}  quantization of
$\big({\rm hsp}(2n,\r),{\rm h}(2n)\big)$ which exponentiates to a
unitary representation of \linebreak[2] $\widetilde{\rm HSp}(2n,\r)$.
\label{prop:unique}
\end{prop}

By ``unique,'' we mean up to unitary equivalence and coextension of
representations (as explained at the end of \S3).

\ms

\noi {\bf Proof:} Suppose $\q$ is a another such quantization of
$\big({\rm hsp}(2n,\r),{\rm h}(2n)\big)$ on some Hilbert space $\h$.
Then
$\q\big({\rm hsp}(2n,\r)\big)$ can be integrated to a representation
$\tau$ of ${\widetilde{\rm HSp}}(2n,\r)$, and (Q4) implies that
$\tau$, when restricted to H$(2n)
\subset {\widetilde{\rm HSp}}(2n,\r)$, is irreducible. The Stone-\vn\
Theorem then states that this representation of H$(2n)$ is unitarily
equivalent to the Schr\"odinger representation, and hence $\tau =
U\varpi U^{-1}$ for some unitary map $U:L^2(\r^n) \ra \h$ by
\cite[Prop. 4.58]{Fo}. Consequently, $\q(f) = U\overline{d\varpi(f)}
U^{-1} | D$ for all $f \in {\rm hsp}(2n,\r)$. Since the Hamiltonian
vector fields of such
$f$ are complete, the corresponding operators
$\q(f)$ are e.s.a., and therefore $\q(f)$ and
$U\overline{d\varpi(f)} U^{-1}$ are coextensive.~$\Box$

\ms

The first, and weakest version of the no-go theorem is:
\begin{thm}[Weak No-Go Theorem] The extended metaplectic quantization
of\/
$(P^2,P^1)$ cannot be extended beyond $P^2$ in $P$.
\label{thm:weak}
\end{thm}

Since $P^2$ is a maximal Poisson subalgebra of $P$ \cite[\S 16]{GS},
we may restate this as: {\em There exists no quantization of\/
$(P,P^1)$ which reduces to the extended metaplectic quantization on
$P^2.$}
\ms

\noi {\bf Proof:} Let $\q$ be a quantization of $(P,P^1)$ which
extends the metaplectic quantization of $(P^2,P^1)$. We will show
that a contradiction arises when cubic polynomials are considered.

By inspection of (\ref{eq:p1}) and (\ref{eq:p2}), we see that the
product
$\ra$ anti-commutator rule (\ref{eq:ac}) is valid for $f,\,g \in P^1.$
In particular, we have the
\vn\ rules
\begin{equation}
\q(q^2)=\q(q)^2,\;\; \q(p^2)=\q(p)^2
\label{eq:vnrr2n1}
\end{equation}

\noi and
\begin{equation}
\q(qp) = {\textstyle \frac{1}{2}}\big(\q(q) \q(p) + \q(p) \q(q)\big).
\label{eq:vnrr2n2}
\end{equation}

\noi These in turn lead to ``higher degree \vn\ rules.''

\begin{lem} For all real-valued polynomials $r$,
\[\q\big(r(q)\big) = r\big(\q(q)\big),\;\;\;\q\big(r(p)\big)=
r\big(\q(p)\big),\]
\[\q\big(r(q)p\big) = {\textstyle
\frac{1}{2}}\big[r\big(\q(q)\big)\q(p)+\q(p)r\big(\q(q)\big)\big],\]

\noi and
\[\q\big(qr(p)\big) = {\textstyle
\frac{1}{2}}\big[\q(q)r\big(\q(p)\big)+r\big(\q(p)\big)\q(q)\big].\]
\label{lem:morevnrs}
\end{lem}

\noi {\bf Proof:} We illustrate this for $r(q) = q^3$. The other rules
follow similarly using induction. Now $\{q^3,q\} = 0$ whence by (Q1)
we have
$\big[\q(q^3),\q(q)\big] = 0$. Since also
$\big[\q(q)^3,\q(q)\big] = 0$, we may write $\q(q^3) =
\q(q)^3 + T$ for some operator $T$ which (weakly) commutes with
$\q(q)$. We likewise have
\[
\big[\q(q^3),\q(p)\big] =  -i\hbar\,
\q\big(\{q^3,p\}\big)
 =  3i\hbar\,\q(q^2) =  3i\hbar\,\q(q)^2 =
\big[\q(q)^3,\q(p)\big]\]

\noi from which we see that $T$ commutes with $\q(p)$ as well.
Consequently, $T$ also commutes with $\q(q)\q(p)+ \q(p)\q(q)$. But
then from (\ref{eq:vnrr2n2}),
\begin{eqnarray*}\q(q^3) & = & {\textstyle \frac{1}{3}}\,
\q\big(\{pq,q^3\}\big) = {\textstyle
\frac{i}{3\hbar}}\,\big[\q(pq),\q(q^3)\big] \\ & =& {\textstyle
\frac{i}{3\hbar}}\,\left[{\textstyle
\frac{1}{2}}\big(\q(q)\q(p) +
\q(p)\q(q)\big),\q(q)^3 + T\right] \\ & = & {\textstyle
\frac{i}{6\hbar}}\,\left[\q(q)\q(p) +
\q(p)\q(q),\q(q)^3\right]
 = \q(q)^3. \;\;\; \bigtriangledown
\end{eqnarray*}

\ms

With this lemma in hand, it is now a simple matter to prove the no-go
theorem. Consider the classical equality

\[{\textstyle \frac{1}{9}}\{q^3,p^3\}= {\textstyle
\frac{1}{3}}\{q^2p,p^2q\}.\]

\noi Quantizing and then simplifying this, the formul{\ae} in
Lemma~\ref{lem:morevnrs} give
\[\q(q)^2\q(p)^2 - 2i\q(q)\q(p) - {\textstyle \frac{2}{3}}I\]

\noi for the L.H.S., and
\[\q(q)^2\q(p)^2 - 2i\q(q)\q(p) - {\textstyle \frac{1}{3}}I\]

\noi for the R.H.S., which is a contradiction. $\Box$

\ms


In \gr's paper \cite{Gr} a much stronger result was claimed; in our
terminology, his assertion was that there is no quantization of
$(P,P^1)$. This is not what Theorem~\ref{thm:weak} states. For if
$\q$ is a quantization of $(P,P^1)$, then while of course $\q(P^1)$
must coincide with Schr\"odinger quantization, $\q$ need {\em not\/}
be the extended metaplectic quantization when restricted to $P^2$.
The problem is that $\q(P^2)$ is not necessarily integrable; (Q5)
only guarantees that
$\q(P^1)$ can be integrated.

With an extra assumption which will guarantee the integrability of
$\q(P^2)$, it is therefore possible to obtain a ``true'' no-go result.
\vh\ supplied such an assumption, which in particular implies: if the
Hamiltonian vector fields of $f,\,g$ are complete and $\{f,g\} = 0$,
then
$\q(f)$ and $\q(g)$ {\em strongly\/} commute. In our view, however,
his assumption is ad hoc, and stronger than what one actually needs
(or wants). We find it preferable to enforce the integrability of
$\q(P^2)$ in a more direct manner. Noting that $P^2$ is the Poisson
normalizer of $P^1$ in $P$, we strengthen (Q5) as follows. Given a
basic set $\bb \subseteq \oo$, let ${\cal N}_{\oo}\big(\wp(\bb)\big)
:= {\cal N}\big(\wp(\bb)\big) \cap \oo$ be the normalizer of
$\wp(\bb)$ {\em in}
$\oo$ \big(where $\nn\big(\wp(\bb)\big)$ denotes the normalizer of
$\wp(\bb)$ in $\p$\big), and set
\begin{enumerate}
\begin{enumerate}
\item[(Q5$'$)] $D$ contains a dense set of separately analytic
vectors for (some Lie generating basis of)
$\q\big({\cal N}_{\oo}(\wp(\bb))\big)$.
\end{enumerate}
\end{enumerate}

\noi We call a quantization satisfying (Q5$'$) in place of (Q5) {\sl
strong}. In view of Remark {\bf 6} (Q5$'$) is not a severe
restriction mathematically, provided ${\cal
N}_{\oo}\big(\wp(\bb)\big)$ is finite-dimensional, and the
alternative -- which ultimately results in a nonintegrable
quantization of ${\cal N}_{\oo}\big(\wp(\bb)\big)$ -- is clearly
pathological. Moreover, as we will see,
${\cal N}_{\oo}\big(\wp(\bb)\big)$ plays an important role in
determining the maximal Poisson subalgebras which are quantizable,
and so it is natural to single it out in this manner. In any case,
with (Q5$'$) we are now able to state a strong no-go result.
\begin{thm}[Strong No-Go Theorem] There does not exist a strong
quantization of
$\big(\p(\r^{2n}),{\rm h}(2n)\big)$.
\label{thm:strong}
\end{thm}

\noi {\bf Proof:} Actually, we will prove a sharper result, viz.,
there exists no strong quantization of $(P,P^1)$. If $\q$ were such a
quantization, then as $P^2$ is finite-dimensional (Q5$'$) implies that
$\q(P^2)$ is integrable and, as before, (Q4) implies that
$\q(P^1)$ is equivalent to the Schr\"odinger representation. Then
Proposition~\ref{prop:unique} shows that
$\q(P^2)$ must be equivalent to the extended metaplectic
quantization, and this contradicts Theorem~\ref{thm:weak}. The
Theorem now follows from this result and the fact that the normalizer
of $P^1$ in $\p$ is the same as its normalizer in $P$, viz. $P^2.$
$\Box$

\ms

Van~Hove \cite{vH1} actually gave a more
refined argument for Theorem~\ref{thm:strong} using in his analysis
only those observables
$f\in \p$ with complete Hamiltonian vector fields, and still obtained
a contradiction from $\q$.

\ms

Finally, we hasten to add that there are  subalgebras of $P$ other
than $P^2$ which can be quantized. For example, let
\[S = \left\{\sum_{i=1}^n f^i(q)p_i + g(q)\right\},\]

\noi where $f^i$ and $g$ are polynomials. Then it is straightforward
to verify that $\sigma:S \ra {\rm Op}(\s (\r^n))$ given by
\begin{equation}
\sigma\big(f(q)p + g(q)\big) = -i\hbar\left(f(q)\frac{d}{dq} +
\frac{1}{2}f'(q)\right) + g(q)
\label{eq:sigma}
\end{equation}

\noi (for $n=1$) is a (strong) quantization of $(S,P^1).$ $\sigma$ is
the familiar ``position'' or ``coordinate representation'' in quantum
mechanics. Since
$S$ is also a maximal subalgebra of $P$, we are able to prove the
following analogue of Theorem~\ref{thm:weak}.
\begin{thm} The position representation $\sigma$ of $(S,P^1)$ cannot
be extended beyond
$S$ in $P$.
\end{thm}

\noi {\bf Proof:} Suppose $\q$ were a quantization of $(P,P^1)$ which
extends
$\sigma.$ (Since $S$ is maximal in $P$, (Q1) implies that any
quantization which extends $\sigma$ must be defined on all of $P$.)
Consider $\q(p^2).$ Mimicking the arguments in the proof of
Lemma~\ref{lem:morevnrs}, we find that $\q(p^2) =
\q(p)^2 + T$, where $T$ commutes with both $\q(q) = q$ and $\q(p) =
-i\hbar\frac{d}{dq}$. Quantizing the classical relation $2p^2 =
\{p^2,pq\}$, we may then show using (\ref{eq:sigma}) that
$T = 0,$ and so
$\q$ restricts to the extended metaplectic representation on $P^2$.
Thus $\q$ must also extend
$d\varpi$ and the result now follows from Theorem~\ref{thm:weak}.
$\Box$

\ms

It unfortunately does not seem possible to prove a uniqueness result
for $\sigma$ analogous to Proposition~\ref{prop:unique}. The reason
stems in part from the fact that ${\cal N}_S(P^1) =
\sp\{1,p,q,q^2,pq\}$ does not satisfy $\big\{{\cal N}_S(P^1),{\cal
N}_S(P^1)\big\} ={\cal N}_S(P^1)$, cf.
\cite[Ch. 4]{Fo}.

 A similar analysis applies to the the Fourier transform of the
subalgebra
$S,$ i.e., the ``momentum'' subalgebra of all polynomials which are at
most affine in the coordinates $q^i.$ In fact, it is not difficult to
see that $P^2$,
$S$ and its Fourier transform exhaust the list of maximal subalgebras
of $P$ which contain $P^1.$

\end{subsection}


\begin{subsection}{$S^2$} Now we turn our attention to the sphere.
Since $S^2$ is compact, all classical observables are complete.
Moreover, the basic set u(2) = span$\{1,\,S_1,\,S_2,\,S_3\}$ is a
compact Lie algebra (cf. \S 5). Consequently all the functional
analytic difficulties present in the case of
$\r^{2n}$ disappear. But the actual computations, which were fairly
routine for
$\r^{2n}$, turn out to be much more complicated for
$S^2.$

The Poisson bracket on $C^{\infty}(S^2)$ corresponding to
(\ref{eq:sfs2}) is
\[\{f,g\}= - \sum_{i,j,k = 1}^3\epsilon_{ijk}\,S_i\,\frac{\partial
f}{\partial S_j} \frac{\partial g}{\partial S_k}.\]

\noi In particular, we have the relations
$\{S_j,S_k\} = -\sum_{l = 1}^3\epsilon_{jkl}\,S_l$. Let $H_k$ denote
the space of spherical harmonics of degree $k$, and define $P^k =
\oplus_{l=0}^k\, H_l$, where the orthogonal direct sum is given by
harmonic decomposition. We may identify elements $f$ of
$P^k$ with polynomials of degree at most $k$ in the components $S_i$
of the spin vector, subject to the relation
\[S_1\,\!^2 + S_2\,\!^2 + S_3\,\!^2 = 1.\]

\noi(By the ``degree'' of such an $f$ we mean the minimum degree of its
representatives.) Note that
$P^1 \cong {\rm u}(2)$. Set
$P = \oplus_{k=0}^{\infty}\,H_k.$

\ms

Let $\q$ be a quantization of $\big(\p(S^2),{\rm u}(2)\big)$ on a
Hilbert space $\h$, whence
\begin{equation}
\big[\q(S_j),\q(S_k)\big] = i\hbar \sum_{l=0}^3
\epsilon_{jkl}\,\q(S_l)
\label{eq:cr}
\end{equation}

\noi and
\begin{equation}
\q({\bf S}^2) = s^2I.
\label{eq:cc}
\end{equation}

\noi By (Q5) and Proposition~\ref{prop:fs}, $\q\big({\rm u}(2)\big)$
can be exponentiated to a unitary representation of the universal
cover SU(2) $\times \r$ of U(2) which, according to (Q4), is
irreducible. Therefore $\h$ must be finite-dimensional, and
$\q\big({\rm u}(2)\big)$  must be one of the usual spin angular
momentum representations, labeled by $j=0,\,\hlf,\,1,\ldots$ For a
fixed value of
$j$, $\dim \h = 2j+1$ and
\begin{equation}
\sum_{i=1}^3\q(S_i)^2=\hbar^2j(j+1)I.
\label{eq:qc}
\end{equation}

Our goal is show that no such (nontrivial) quantization exists.
Patterning our analysis after that for $\r^{2n}$, we use
irreducibility to derive some generalized
\vn\ rules.
\begin{lem} For $i = 1,\,2,\,3$ we have
\begin{equation}
\q\big(S_i\,\!^2\big)=a\q(S_i)^2+cI
\label{eq:sisi}
\end{equation}

\noi where $a,\,c$ are representation dependent real constants with
$ac \neq 0$.
\label{lem:vnrs21}
\end{lem}

The proof is in \cite{GGH}. {}From this we also derive
\begin{equation}
\q(S_iS_k)={\displaystyle
\frac{a}{2}}\big(\q(S_i)\q(S_k)+\q(S_k)\q(S_i)\big)
\label{eq:sisk}
\end{equation}

\noi for $i \neq k$. (As an aside, these formul{\ae} show that a
quantization, if it exists, is badly behaved with respect to the
multiplicative structure on
$C^{\infty}(S^2)$; in particular, the product $\ra$ anti-commutator
rule cannot hold. Remarkably, this is as it should be: for {\em if\/}
this rule were valid, then -- subject to a few mild assumptions on
$\q$ -- the classical spectrum of $S_3$, say, would have to coincide
with that of $\q(S_3)$ which is contrary to experiment \cite{GGH}.)
With these tools, we can now prove the main result:
\begin{thm} There is no nontrivial quantization of $\big(\p(S^2),{\rm
u}(2)\big)$.
\label{thm:nogos2}
\end{thm}

\noi {\bf Proof:}  As for Theorem~\ref{thm:strong}, we will actually
prove a sharper result, viz. there is no nontrivial quantization of
$(P,P^1)$.

Fix $j>0$, as $j=0$ produces a trivial quantization. Assuming that
$\q$ is a quantization of $(P,P^1)$, we can use
(\ref{eq:cr})-(\ref{eq:sisk}) to quantize the classical relation
\[s^2S_3=\big\{S_1\,\!^2-S_2\,\!^2,\,S_1S_2\big\}-\big\{S_2S_3,\,S_3S_1\big\},\]

\noi thereby obtaining
\begin{equation} s^2=a^2\hbar^2\big(j(j+1)-{\textstyle
\frac{3}{4}}\big)
\label{eq:xxx}
\end{equation}

\noi which  contradicts $s>0$ for $j=\hlf$. Now assume $j>\hlf$, and
quantize
\[ 2s^2S_2S_3=\big\{S_2\,\!^2,\{S_1S_2,S_1S_3\}\big\} - \f 3,4.
\big\{S_1\,\!^2,\{S_1\,\!^2,S_2S_3\}\big\},\]

\noi similarly obtaining
\[s^2=a^2\hbar^2\big(j(j+1)-\f 9,4.\big)\]

\noi which contradicts (\ref{eq:xxx}). Thus we have derived
contradictions for all $j>0$, and the theorem is proven.
$\Box$

\ms

In view of the impossibility of quantizing $(P,P^1)$, one can ask
what the maximal subalgebras in $P$ are to which we can extend an
irreducible representation of $P^1$. The following chain of results,
which we quote without proof (cf. \cite{GGH}), provides the answer.
\begin{prop}
$P^1$ is a maximal Poisson subalgebra of
$O \oplus\r\subset P$, where $O$ is the Poisson algebra consisting of
polynomials containing only terms of odd degree.
\label{prop:max}
\end{prop}

Next we establish a no-go theorem for $(O \oplus \r,P^1)$. However,
the generalized Von~Neumann rules listed in Lemma~\ref{lem:vnrs21}
involve only even degree polynomials, so these are not applicable in
$O$. Fortunately, we have another set of generalized Von~Neumann
rules, also implied by the irreducibility of
$\q(P^1)$, involving only terms of odd degree.
\begin{lem} If $\q$ is a quantization of $(O \oplus\r,P^1)$, then for
$i=1,\,2,\,3$,
\[\q\big(S_i\,\!^3\big)=b\q(S_i)^3+e\q(S_i)\]

\noi where $b,\,e\in\r$.
\end{lem}

{}From this we prove (with far greater effort):
\begin{thm} There is no nontrivial quantization of
$(O\oplus\r,P^1)$.
\label{thm:yyy}
\end{thm}

Now $O \oplus \r$ is itself a maximal subalgebra of $P$, and in fact
the only Poisson subalgebras of $P$ strictly containing $P^1$ are
$O \oplus \r$ and $P$ itself. Thus Theorem~\ref{thm:yyy} and
Proposition~\ref{prop:max} combine to yield our sharpest result for
the sphere:
\begin{cor} No nontrivial quantization of $P^1$ can be extended to a
larger subalgebra of $P$.
\end{cor}

An interesting observation is that  there is {\em no} obstruction for
the quantum spin $j=0$ representation. In fact, there exists a
(unique) trivial quantization of $(\p,P^1)$ with
$\q(S_i)=0$, but $\q\big(S_i\,\!^2\big)=\frac{s^2}{3} I$ for all $i$. It is
defined by
$Q(f)=f_0$, where $f_0$ is the constant term in the harmonic
decomposition of $f \in \p$.

There are crucial structural differences between the
Groenewold-Van~Hove analysis of $\r^2$ and the current analysis of
$S^2$. Within $P$, $\sp\{1,\,p,\,q\}$ has as its Poisson normalizer
the algebra of polynomials at most degree 2, and there is no
obstruction to quantization in this algebra: the obstruction comes
from the cubic polynomials. On the other hand, for the sphere, the
algebra $\sp\{1,\,S_1,\,S_2,\,S_3\}$ is self-normalizing in $P$; we
obtain an obstruction in the quadratic polynomials, and find that
there is no extension possible for a quantization of $P^1$. The fact
that this u(2)-subalgebra is self-normalizing is one reason why we
are able to obtain ``strong'' no-go results for the sphere relatively
easily (as compared to $\r^{2n}$).

\end{subsection}


\begin{subsection}{$T^2$}

At the end of the previous subsection, we showed that there was a full
quantization of the sphere, albeit a ``trivial'' one. Here we exhibit
a nontrivial full quantization of the Poisson algebra of the torus.
Proofs for the results in this section can be found in \cite{Go}.

Consider the torus $T^2$ thought of as $\r^2/\z^2$, with symplectic
form
$$\omega=B\,dx\wedge dy.$$ We will study the family
\[{\cal B}_k = \mbox{\rm span}\{1,\sin 2\pi kx,\cos 2\pi kx,\sin 2\pi
ky,\cos 2\pi ky\}\]

\noi of basic sets with $k$ a positive integer. The crucial difference
between this example and the previous ones is that the Poisson
algebras
$\wp(\bb_k)$ are all infinite-dimensional.

Now $(T^2,\omega)$ is quantizable provided $B = N\hbar$ for some
nonzero integer $N$. Fix $N = 1$ and let $L$ be the corresponding
Kostant-Souriau prequantization line bundle over
$T^2$ \cite{ki}. Then one can identify the space of smooth sections
$\Gamma(L)$ with the space of ``quasi-periodic'' functions $\varphi\in
C^\infty(\r^2)$ satisfying
\[\varphi(x+m,\,y+n)=e^{2\pi imy}\varphi(x,\, y)\,,\quad n,\,
m\in\z,\]

\noi and the prequantization Hilbert space $\h$ with the (completion
of) the set of those quasi-periodic $\varphi$  which are
$L^2$ on
$[0,\,1)\times[0,\, 1)$. The associated prequantization map
$\q:\p\to{\rm Op}\big(\Gamma(L)\big)$ (for a specific choice of
connection on $L$) is defined by
\[\q(f)=-i\hbar\left[{\partial f\over\partial x}\bigg(
{\partial\over\partial y}-\frac{i}{\hbar}x\bigg)-{\partial
f\over\partial y} {\partial\over\partial x}\right]+f.\]

\noi As the torus is compact, these operators are essentially
self-adjoint on $\Gamma(L)\subset\h$.

\begin{thm}
$\q$ is a quantization of $(\p,\bb_k)$ for all positive integers $k$.
\label{thm:nogot2}
\end{thm}

\noi {\bf Proof:}  Since $\q$ is a prequantization, it suffices to
verify (Q4) and (Q5). To this end it is convenient to use complex
notation and write
\[\bb_k = \sp \!\!\left\{1,e^{\pm 2\pi ikx},e^{\pm 2\pi
iky}\right\}.\]

The analysis is simplified by applying the Weil-Brezin-Zak transform
$Z$ \cite[\S1.10]{Fo} to the above data. Define a unitary map $Z: \h
\rightarrow L^2(\r)$ by
\[(Z\phi)(x) = \int_0^1\phi(x,y)\,dy\]

\noi with inverse
\[(Z^{-1}\psi)(x,y) = \sum_{m \in {\bf Z}}\psi(x+m)e^{-2\pi imy}.\]

\noi Under $Z$ the domain $\Gamma(L)$ maps onto the Schwartz space $\s
(\r)$ \cite{ki}. Setting $A_{\pm} :=Z\q(e^{\pm 2\pi ikx})Z^{-1}$ and
$B_{\pm} :=Z\q(e^{\pm 2\pi iky})Z^{-1}$ we compute, as operators on
$\s (\r)$,
\bea (A_{\pm}\psi)(x)  & = & e^{\pm 2\pi ikx}(1 \mp 2\pi ikx)\psi(x)
\nonumber
\\ (B_{\pm}\psi)(x) & = & \bigg(1 \mp 2 \pi \hbar k
\frac{d}{dx}\bigg)\psi(x \pm k). \nonumber
\eea

\noi Then
$A_{\pm}\,\!^* = \ol{A_{\mp}}$ on the domain $\{\psi\,|\, x\psi \in
L^2(\r)\}$, and likewise $B_{\pm}\,\!^* =
\ol{B_{\mp}}$ on $\{\psi\,|\, d\psi/dx \in
L^2(\r)\}$.\footnote{\,$d\psi/dx$ is to be understood in the sense of
tempered distributions.} In fact
$\ol {A_{\pm}}$ and
$\ol {B_{\pm}}$ are normal operators.

To show that $\q(\bb_k)$ is an irreducible set, let us suppose that
$T$ is a bounded s.a. operator on $L^2(\bf R)$ which strongly commutes
with $\ol {A_{\pm}}$ and $\ol {B_{\pm}}$. Then $T$ must commute (in
the weak sense) with these operators on their respective
domains.\footnote{\,Here and in what follows we use the fact that a
bounded operator weakly commutes with an (unbounded) normal operator
iff they strongly commute.} Consequently $T$ commutes with both
\[{\ol{A_-}}\,{\ol{A_+}} = I + 4 \pi^2 k^2 x^2\]

\noi  on the domain $\{\psi\,|\,x^2\psi \in L^2(\r)\}$, and

\[{\ol{B_-}}\,{\ol{B_+}} = I - 4\pi^2 \hbar^2 k^2 \frac{d^2}{dx^2}\]

\noi on $\{\psi\,|\,d^2\psi/dx^2 \in L^2(\r)\}$. {}From these
equations we see that
$T$ commutes, and therefore strongly commutes, with the closures of
two of the three generators of the metaplectic representation
(\ref{eq:p2}) of sp$(2,\r)$ on
$\s (\r)$.

Let $\mu$ denote the metaplectic representation of the metaplectic
group \linebreak[4] Mp(2,$\r$) on
$L^2(\r)$. We have in effect just established that $T$ commutes with
the one parameter groups
$\exp\!\big(is\,{\ol{x^2}}\big)$ and
$\exp\!\big(\!-\!it\hbar^2\,{\ol{d^2/dx^2}}\big)$. Now classically
the exponentials
$\exp(sx^2)$ and
$\exp(ty^2)$ generate Sp$(2,\r)$ \cite[\S 4]{GS}. As Mp$(2,\r) \ra
{\rm Sp}(2,\r)$ is a double covering, the corresponding exponentials
in Mp(2,$\r$) generate a neighborhood of the identity in the
metaplectic group.  Since $\mu\big[\exp(sx^2)\big] =
\exp\!\big(is\,{\ol{x^2}}\big)$ and
$\mu\big[\exp(ty^2)\big] =
\exp\!\big(\!-\!it\hbar^2\,{\ol{d^2/dx^2}}\big)$, it follows that $T$
commutes with $\mu({\cal M})$ for all ${\cal M}$ in a neighborhood of
the identity in Mp(2,$\r$) and hence, as this group is connected, for
all
${\cal M} \in$ Mp$(2,{\bf R})$.

Although the metaplectic representation $\mu$ is reducible, the
subrepresentations
$\mu_e$ and $\mu_o$ on each invariant summand of $L^2({\bf R})  =
L^2_{e}({\bf R})
\oplus L^2_{o}({\bf R})$ of even and odd functions are irreducible
\cite[\S 4.4]{Fo}. Writing $T = P_eT + P_oT$, where $P_e$ and $P_o$
are the even and odd projectors, one has
\begin{equation} [P_eT,\mu({\cal M})] = 0 \label{eq:com}
\end{equation}

\noi for any ${\cal M} \in \mbox{Mp(2,}{\bf R})$. It then follows
from the irreducibility of the subrepresentation $\mu_e$ that
$P_eT = k_eP_e + RP_o$ for some constant $k_e$ and some operator
$R:L^2_{o}({\bf R}) \ra L^2_{e}({\bf R})$. Substituting this
expression into (\ref{eq:com}) yields $[RP_o,\mu({\cal M})] = 0$, and
Schur's Lemma then implies that $R$ is either an isomorphism or is
zero. But $R$ cannot be an isomorphism as the representations $\mu_e$
and $\mu_o$ are inequivalent \cite[Thm.~4.56]{Fo}. (Recall that two
unitary representations are similar iff they are unitarily
equivalent.) Thus
$P_eT = k_eP_e$. Similarly
$P_oT = k_oP_o,$ whence $T = k_eP_e + k_oP_o.$

But now a short calculation shows that $T$ commutes with

\[{\ol {A_+}}-{\ol{A_-}} = 2i(\sin 2\pi kx - 2\pi kx\cos 2\pi kx)\]

\noi only if $k_e = k_o$. Thus $T$ is a multiple of the identity, and
so $\{A_{\pm},B_{\pm}\}$ is an irreducible set, as was to be shown.
Thus (Q4) is satisfied.

For (Q5), we claim that the linear span of the Hermite functions form
a dense set of separately analytic vectors for the e.s.a.\ components
of
$\{A_{\pm},B_{\pm}\}$. {}From the expression above for $A_{\pm}$, it
is clear that a vector will be analytic for the e.s.a.\ components of
$A_{\pm}$ iff it is analytic for multiplication by $x$. But it is
well known that the Hermite functions are analytic for this latter
operator. The result for
$B_{\pm}$ is obtained directly from this by means of the Fourier
transform.{}~$\Box$

\ms

\noi {\bf Remark: 10.} It is interesting to note that the proof breaks
down when
$|N| \neq 1$ \cite{Go}. While it is not known to what extent this
theorem will remain valid in general (but see \S 5), one can prove in
the special case when $k = |N|$ that the corresponding
prequantization map does {\em not} represent $\bb_N$ irreducibly, and
so does not qualify as a quantization of
$\big(\p(T^2,\omega_N),\bb_N\big)$.

\end{subsection}

\end{section}


\begin{section}{Speculations}

Comparing the behavior of the examples presented in the previous
section, we see that we were able to obtain a very sharp no-go
theorem for
$S^2$, a relatively weaker no-go result for $\r^{2n}$,  and no
obstruction at all in the case of $T^2$. On the basis of these
examples, we attempt to extract the key features which govern the
appearance of obstructions to a full quantization as well as their
severity.

Of course, any conclusions that we can draw at this point are
necessarily quite tentative, due to the paucity of examples against
which to test them. There are also various aspects of our three main
examples that are still not well understood. Nonetheless, some
interesting observations can be made, which may prove helpful in
subsequent investigations.

Of our three examples, the torus is clearly much different than the
others. It is not a Hamiltonian homogeneous space, and the Poisson
algebras generated by the basic sets $\bb_k$ are
infinite-dimensional. The latter property seems to be the controlling
factor: In effect, since $\wp(\bb_k)$ is infinite-dimensional, the
irreducibility requirement (Q4) loses much of its force -- so much
so, in the case of the torus, that it precludes the existence of an
obstruction. Hence we propose that a general Groenewald-Van~Hove
theorem takes the form:
\begin{conj} Let $M$ be a symplectic manifold and
$\bb  \subset \p(M)$ a basic set with $\wp(\bb)$ finite-dimensional.
Then there is no nontrivial strong quantization of \linebreak[2]
$\big(\p(M),\bb\big)$.
\label{conj:gvh}
\end{conj}

A distinction must be made here depending upon whether $\wp(\bb)$ is a
``compact'' Lie algebra (i.e., is the Lie algebra of a compact Lie
group)\footnote{\, This is nonstandard terminology.} or not. The
reason is that by (Q4) and (Q5) the representation
$\q\big(\wp(\bb)\big)$ will be integrable to an irreducible unitary
representation ${\mit \Pi}$ of the connected, simply connected Lie
group
$\tilde G$ whose Lie algebra is $\wp(\bb)$. If $\wp(\bb)$ is compact,
then some discrete quotient $G$ of $\tilde G$ will be compact. Let
$\Gamma$ be the kernel of the covering map $\tilde G \ra G$; then
$\Gamma$ is central in
$\tilde G$. By irreducibility, ${\mit \Pi}(\gamma)$ must be a
multiple of the identity for each $\gamma \in \Gamma$ and so, by
unitarity, each ${\mit
\Pi}(\gamma)$ has unit modulus.  It follows that ${\mit \Pi}$  can be
viewed as a projective representation of $G.$ This projective
representation, in turn, can be realized as an honest representation
of a certain central extension of
$G$ by ${\mit \Pi}(\Gamma)$ \cite[\S 13.2]{b-r}. Since ${\mit
\Pi}(\Gamma) \subset T^1$ and is discrete, it is finite, and hence
this central extension is compact. Consequently the irreducible
representations of the latter are finite-dimensional, and so all
functional analytic difficulties are obviated. In the compact case,
then, one can delete the adjective ``strong'' in
Conjecture~\ref{conj:gvh}. This explains (in part) why it was
``easier'' to obtain the no-go result for the sphere than for
$\r^{2n}$.

In particular, generalizing the analysis of $S^2$ given in
\cite{GGH}, it appears likely that ${\bf C}P^n$ with basic set u$(n)$
will admit no nontrivial full quantizations.

When $\wp(\bb)$ is noncompact it is necessary to confront these
analytic difficulties. In fact, it is the noncompactness of h($2n$)
that is responsible for the possible non-integrability of $\q(P^2)$,
and hence the splitting of the no-go results for $\r^{2n}$ into
``weak'' and ``strong'' cases. In this regard, it would be very
interesting to determine whether there exists a full {\em weak\/}
quantization of
$\big(\p(\r^{2n}),{\rm h}(2n)\big)$, or at least a weak quantization
of
$(P,P^1)$. Of course, $\q(P^2)$ cannot then be integrable, and the
\vn\ rules (\ref{eq:vnrr2n1}) cannot hold. The physical import of
this is discussed in \cite{KLZ}.

As all extant no-go results are for elementary systems, it would be
useful to study symplectic manifolds which are not Hamiltonian
homogeneous spaces, but which nonetheless have basic sets with
$\wp(\bb)$ finite-dimensional. A particularly interesting example of
this type is
$T^*S^1$, with
$\bb = \sp\{1,\sin \theta,\cos \theta,l\}$, where $(\theta,l)$ are
canonical coordinates with $l$ being the angular momentum. It lies
``halfway'' between $\r^2$ and
$T^2$, and so should provide a good test of Conjecture~\ref{conj:gvh}.
Other examples worth studying  are various elementary systems for
Sp($2n,\r) \times \r$, such as
$\r^{2n}\setminus \{\bf 0\}$ and its coadjoint orbits.\footnote{\, The
factor $\r$ is inserted here so that the corresponding basic set will
include the constants, cf. (B4).}

In view of the torus example it seems equally reasonable to propose
\begin{conj} Let $M$ be a symplectic manifold and
$\bb$ a basic set with $\wp(\bb)$ dense in $C^{\infty}(M)$. Then
there exists a nontrivial quantization of $\big(\p(M),\bb\big)$.
\label{conj:go}
\end{conj}

Obviously, a necessary condition for $\q$ to be a full quantization
of $(\p,\bb)$ is that $\q$ represent $\p$ itself irreducibly. It
turns out \cite{ch3,Tu95} that this is so for all Kostant-Souriau
prequantizations\footnote{\,However, it should be noted that there
are other prequantizations which do not represent $\p$ irreducibly,
for instance, the prequantization of Avez \cite{av,Ch2}.}; thus it is
natural to consider the case when $M$ is prequantizable in this sense.
In fact, in this context \cite{Tu95} gives even more:
\begin{prop} Let $M\!$ be an integral symplectic manifold, $\!L$ a
Kostant-Souriau prequantization line bundle over $M$ and
$\q_L$ the corresponding prequantization map. Let
$\bb$ be a basic set with
$\wp(\bb)$ dense in $C^{\infty}(M)$. Then $\q_L$ represents $\bb$
irreducibly on the domain which consists of compactly supported
sections of $L$.
\label{prop:gijs}
\end{prop}

Set $D_0 =
\Gamma(L)_0$, the compactly supported sections of $L$. By
construction $\q_L:\p(M) \ra {\rm Op}(D_0)$ satisfies (Q1)-(Q3). This
proposition states that $\q_L$ satisfies (Q4) as well. Thus to obtain
a full quantization it remains to verify (Q5) -- perhaps on some
appropriately chosen coextensive domain $D$; unfortunately, it does
not seem possible to do this except in specific instances. A first
test would be to understand what happens for
$\big(\p(T^2,\omega_N),\bb_1\big)$. In any event,
Proposition~\ref{prop:gijs} does provide a certain amount of support
for Conjecture~\ref{conj:go}.

The ``gray area'' between these two conjectures consists of
symplectic manifolds with basic sets $\bb$ for which $\wp(\bb)$ is
infinite-dimensional, yet not dense in
$C^{\infty}(M)$. That the infinite-dimensionality of
$\wp(\bb)$ alone may be enough to guarantee the existence of a full
quantization is evinced by our results for
$\big(\p(T^2,\omega_1),\bb_k\big)$ with $k>1.$ (See however Remark
{\bf 10}.)

Proving Conjecture~\ref{conj:gvh} in any sort of generality seems
well beyond what is possible with present technology. Still, based on
the examples that have already been worked out, we can suggest an
avenue of attack: Assume that there does exist a full quantization of
$(\p(M),\bb)$. The irreducibility and integrability  requirements
upon $\q(\bb)$ should lead to \vn\ rules for elements of $\bb$. For
example, in the case of $\r^{2n}$ these rules are:

\[\q\big(x^2\big) = \q(x)^2,\;\;x\in {\rm h}(2n).\]

\noi Similarly, for the sphere, we computed the somewhat more
complicated rules

\[\q\big(S_i\,\!^2\big) = a\q(S_i)^2 + cI,\;\;i = 1,\,2,\,3\]

\noi with $ac \neq 0.$ The idea is then to show that eventually, the
\vn\ rules lead to a contradiction with the Poisson bracket
$\rightarrow$ commutator rule. The main difficulty with this approach
is that the derivation of these rules requires a detailed knowledge
of the representation theory of $\wp(\bb)$. Perhaps it is possible to
streamline this process, at least for certain types of basic sets.
Another problem is that one doesn't know where to look for the
contradiction. For $\r^{2n}$, it arises when one considers cubic
polynomials, while for $S^2$ the contradiction comes sooner, when one
considers spherical harmonics of degree 2. Fortunately, another
conjecture (see below) enables one to ``guess'' where the
contradiction might lie. It will likely prove necessary to work
through a few more examples of
\gr-\vh\ obstructions before one is able to refine this approach
enough to become workable.

As well, it would be useful to consider basic sets other than the
ones we have studied in the standard examples, for instance, the
Euclidean algebras that arise in optics \cite[\S17]{GS}. We have also
restricted consideration to polynomial subalgebras to a large extent,
but there are other subalgebras $\oo$ which are of interest (e.g., on
$\r^{2n}$, those functions which are constant outside some compact
set \cite{Ch2}).

A negative answer to the conjecture might indicate that one should
strength\-en the conditions defining a basic set by, e.g., replacing
(B3) by (C2) as discussed in \S 3. Or, if the conjecture seems
undecidable, perhaps one should abandon the definition of a
quantization map solely in terms of basic sets and consider an
alternative. However, the two other ways to define a quantization map
listed previously suffer from serious flaws. If one imposes
\vn\ rules at the outset, then one tends to run into difficulties
rather quickly -- especially if one tries to enforce the rules on all
of $\p(M)$ and not some basic subset thereof -- as was shown in
\S 4.1. Furthermore, it is unclear what form \vn\ rules should take
in general, as is illustrated by the unintuitive rules for the sphere
exhibited above. For instance, mimicking the situation for
$\r^{2n},$ one might simply postulate that
$\q(x^2) =
\q(x)^2$ for $x \in u(2).$\footnote{\, In effect, this is exactly
what transpired in the case of $\r^{2n}$: condition (Q5$'$) was
specifically designed so as to produce the squaring rule for elements
of h(2$n$).} But the squaring rule for angular momentum is plainly
incompatible with the (presumably) correct
\vn\ rules (\ref{eq:sisi}), because of the restriction that $ac \neq
0.$\footnote{\,Because of this, \cite{KLZ} would refer to
(\ref{eq:sisi}) as ``{\em non\/}-Neumann rules''! } Even disregarding
this point, one would still ``miss'' various possibilities
(corresponding to the freedom in the choice of parameters
$a,\,c$), which do occur in specific representations.  And in the
case of the torus,
\vn\ rules are effectively moot, since the explicit prequantization
map
$\q$ itself determines the quantization of every observable. Another
way of phrasing this is that the irreducibility requirement is so
weak in this context that it cannot enforce any \vn\ rules, and hence
does not lead to a contradiction with the given
$\q.$ All in all, it appears as if the \vn\ rules play a secondary
role; the basic set $\bb$ is the primary object. It is also more
compelling physically and pleasing {\ae}sthetically to require $\q$
to satisfy an irreducibility requirement than a \vn\ rule.

There are problems with the polarization approach as well. For one
thing, symplectic manifolds need not be polarizable \cite{go87}. This
relatively rare occurrence notwithstanding, there are quantizations
which cannot be obtained by polarizing a prequantization: a well
known example is the extended metaplectic quantization of the pair
\big(hsp(2$n$,\r),h(2$n$)\big)
\cite{Bl}. As we shall see presently, the specific predictions of
geometric quantization theory are also off the mark in a number of
instances.

Finally, it should be mentioned that these three approaches to
quantization typically lead to obstructions in one way or another. We
have already seen in the beginning of \S 4.1 that the imposition of
\vn\ rules, in and by themselves, results in inconsistencies;
moreover, \vn\ rules played a crucial role in deriving the \gr-\vh\
obstructions for
$\r^{2n}$ and $S^2$. In the context of polarizations, it is known that
the only observables which are consistently quantizable {\it a
priori} are those whose Hamiltonian vector fields preserve a given
polarization
\cite{bl1,Wo}. While this does not preclude the possibility of
quantizing more general observables, attempts to quantize observables
outside this class in specific examples usually result in
inconsistencies. In {\it all\/} instances, the set of {\it a priori\/}
quantizable observables relative to a given polarization forms a
proper subalgebra of the Poisson algebra of the given symplectic
manifold. This observation provides further corroboration that
\gr-\vh\ obstructions to quantization should be ubiquitous.

\ms

Setting aside the question of the existence of \gr-\vh\ obstructions,
let us now suppose that there is such an obstruction, so that it is
impossible to consistently quantize all of $\p$. The question is:
what are the maximal subalgebras $\oo \subset \p$ containing the
given basic set $\bb$ such that
$(\oo,\bb)$ can be quantized? Now, technical issues aside, given a
representation
$\q$ of
$\wp(\bb)$ on a Hilbert space $\h$, one ought to be able to induce a
representation of its Poisson normalizer ${\nn}\big(\wp(\bb)\big)$ on
$\h$. (Indeed, the structure $\big(\nn(\wp(\bb)),\bb\big)$ brings to
mind an infinitesimal version of a Mackey system of imprimitivity
\cite{b-r}.) In particular, for
$\r^{2n}$ one has that
$\nn\big({\rm h}(2n)\big) = {\rm hsp}(2n,\r)$, and for $S^2$ one
computes that
$\nn\big({\rm u}(2)\big) = {\rm u}(2)$. In both cases, we have shown
that these normalizers are in fact the maximal subalgebras that can be
consistently quantized. Thus it seems reasonable to assert:
\begin{conj} Let $\bb$ be a basic set, with $\wp(\bb)$
finite-dimensional. Then every integrable irreducible representation
of
$\wp(\bb)$ can be extended to a quantization of
$\big(\nn(\wp(\bb)),\bb\big)$. Furthermore, no nontrivial
quantization of
$\big(\nn(\wp(\bb)),\bb\big)$  can be extended beyond
$\nn\big(\wp(\bb)\big)$.
\label{conj:obs}
\end{conj}

If true, this conjecture would point where to look for a \gr-\vh\
contradiction, viz. just outside the normalizer. This is in exact
agreement with the examples. However, it should be emphasized that
there may exist subalgebras $\oo \subset \p$, with
$(\oo,\bb)$ maximally quantizable, which do not arise in this
fashion. A familiar example is the position subalgebra $S =
\big\{f(q)p + g(q)\big\}$ of
$\p(\r^{2})$ encountered in \S 4.1.  It is not clear how one could
``discover'' this subalgebra given just the basic set $\bb = {\rm
h}(2)$ (but see below); note also that
$\nn\big({\rm h}(2)\big) \not\subset S$.

This is reminiscent of the situation in geometric quantization with
respect to polarizations. Suppose that $\cal A$ is a polarization of
$\p_{\bf C}(M)$. Then one knows that one can consistently quantize
those observables which preserve $\cal A$, i.e., which belong to the
real part of $\nn(\cal A)$ \cite{bl1,Wo}. In this way one obtains a
``lower bound'' on the set of quantizable functions for a given
polarization. If one takes the antiholomorphic polarization on $S^2$,
then it turns out that the set of {\em a priori} quantizable
functions obtained in this manner is precisely the u(2) subalgebra
$\sp\{1,\,S_1,\,S_2,\,S_3\}$. But it may happen that the real part of
$\nn(\cal A)$ is too small, as happens for $\r^{2n}$ with the
antiholomorphic polarization. In this case the real part of $\nn(\cal
A)$ is only a proper subalgebra of $P^2$, and in particular is not
maximal. This illustrates the fact, alluded to previously, that the
extended metaplectic representation cannot be derived via geometric
quantization. Furthermore, in the case of the torus, introducing a
polarization will drastically cut down the set of {\em a priori}
quantizable functions, which is at odds with the existence of a full
quantization of this space. So geometric quantization is not a
reliable guide insofar as computing maximally quantizable subalgebras
of observables. On the other hand, the position subalgebra $\s =
\{f(q)p + g(q)\}$ is just the normalizer of the vertical polarization
${\cal A} = \{h(q)\}$ on
$\r^{2n}$, so this subalgebra finds a natural interpretation in the
context of polarizations.

Clearly, there must be some connection between polarizations and
basic sets that awaits elucidation. It would be interesting to
determine if there is a way to recast the \gr-\vh\ results in terms of
polarizations. It would also be worthwhile, assuming
Conjecture~\ref{conj:obs} to be correct,  to see whether one can use the
{\em ab initio} knowledge of a maximal set of quantizable observables to
refine geometric quantization theory, or to develop a new quantization
procedure, which is adapted to the \gr-\vh\ obstruction in that it
will automatically be able to quantize this maximal set.

\ms

One of our goals in this paper was to obtain results which are
independent of the particular quantization scheme employed, as long
as it is Hilbert-space based. Therefore it is interesting that the go
and no-go results described in this proposal have direct analogues in
deformation quantization theory, since this theory
was developed, at least in part, to avoid the use of Hilbert spaces
altogether \cite{Bayen e.a.}. So for example,  the no-go result for
$S^2$ is mirrored by the fact that there are no strict
SU(2)-invariant deformation quantizations  of
$\p(S^2)$
\cite{Ri}, while the go theorem for $T^2$ has as a counterpart the
result that there do exist strict deformation quantizations of the
torus \cite{Ri}. It is generally believed that the existence of
\gr-\vh\ obstructions necessitates a weakening of the Poisson bracket
$\ra$ commutator rule (by insisting that it hold only to order
$\hbar$), but the analogies above indicate that this may not suffice
to remove the obstructions. There are undoubtedly many important
things to be learned by getting to the heart of this phenomenon.

\end{section}



\end{document}